\documentclass[10pt]{amsart}
\usepackage{amssymb,amsmath,amsthm}
\usepackage{mathrsfs,a4wide}
\theoremstyle{plain}
\newtheorem{theorem}{Theorem}[section]
\newtheorem{lemma}[theorem]{Lemma}
\newtheorem{proposition}[theorem]{Proposition}

\newtheorem{remark}[theorem]{Remark}

\theoremstyle{definition}
\theoremstyle{remark}
\numberwithin{equation}{section}

\newcommand{\hs}{{\mathcal H}}

\newcommand{\fs}{{\mathcal F}}

\newcommand{\ns}{{\mathcal N}}
\newcommand{\bs}{{\mathcal B}}

\newcommand{\Es}{{\mathcal E}}

\newcommand{\R}{{\mathbb R}}
\newcommand{\N}{{\mathbb N}}

\newcommand{\tint}[1]{{\rm int}(#1)}

\newcommand{\Om}{\Omega}
\newcommand{\Omb}{\overline{\Omega}}

\newcommand{\weak}{\rightharpoonup}

\newcommand{\mmeps}[1]{\frac{\varepsilon}{2}\int_{\Om} |\nabla#1|^2 \,dx+\frac{1}{2\varepsilon}\int_{\Om}(1-#1)^2 \,dx}

\newcommand{\mmn}[1]{\frac{\varepsilon_n}{2}\int_{\Om} |\nabla#1|^2 \,dx+\frac{1}{2\varepsilon_n}\int_{\Om}(1-#1)^2 \,dx}
\newcommand{\mmnp}[1]{\frac{\varepsilon_n}{2}\int_{\Om'} |\nabla#1|^2 \,dx+\frac{1}{2\varepsilon_n}\int_{\Om'}(1-#1)^2 \,dx}
\newcommand{\elleps}[2]{\int_{\Om} (\eta_\varepsilon+{#2}^2)|\nabla#1|^2 \,dx}
\newcommand{\ellepsn}[2]{\int_{\Om} (\eta_{\varepsilon_n}+{#2}^2)|\nabla#1|^2 \,dx}
\newcommand{\ellepsnp}[2]{\int_{\Om'} (\eta_{\varepsilon_n}+{#2}^2)|\nabla#1|^2 \,dx}
\newcommand{\elln}[2]{\int_{\Om} (\eta_n+{#2}^2)|\nabla#1|^2\,dx}
\newcommand{\ellnp}[2]{\int_{\Om'} (\eta_n+{#2}^2)|\nabla#1|^2\,dx}
\newcommand{\ateps}[2]{\int_{\Om} (\eta_\varepsilon+{#2}^2)|\nabla#1|^2 \,dx+\mmeps{#2}}

\newcommand{\Ateps}[2]{F_\varepsilon(#1,#2)}
\newcommand{\Atepsn}[2]{F_{\varepsilon_n}(#1,#2)}
\newcommand{\Atn}[2]{F_{\varepsilon_n}(#1,#2)}
\newcommand{\ueps}{u_\varepsilon}
\newcommand{\veps}{v_\varepsilon}
\newcommand{\uepsgn}{u_{\varepsilon,h}}
\newcommand{\vepsgn}{v_{\varepsilon,h}}
\newcommand{\uepsn}{u_{\varepsilon_n}}
\newcommand{\vepsn}{v_{\varepsilon_n}}
\newcommand{\ind}[1]{1_{\{#1\}}}

\newcommand{\supt}[1]{{\rm supt}(#1)}

\title
[AMBROSIO-TORTORELLI APPROXIMATION OF QUASI-STATIC EVOLUTION]
{AMBROSIO-TORTORELLI APPROXIMATION OF \\ QUASI-STATIC EVOLUTION OF BRITTLE FRACTURES}
\author[A. Giacomini]
{Alessandro Giacomini}
\address[Alessandro Giacomini]{S.I.S.S.A., Via Beirut 2-4, 34014, Trieste,
Italy}
\email[A. Giacomini]{giacomin@sissa.it}
\begin{document}
\vskip .2truecm
\begin{abstract}
\small{We define a notion of quasi-static evolution for the elliptic approximation of the
Mumford-Shah functional proposed by Ambrosio and Tortorelli. Then we prove that this regular
evolution converges to a quasi static growth of brittle fractures in linearly elastic bodies.
}
\end{abstract}
\maketitle
{\small \tableofcontents}

\section{Introduction}
\label{intr}

In 1998 Francfort and Marigo \cite{FM} proposed a model of quasi-static growth of brittle
fractures in linearly elastic bodies based on the classical Griffith criterion.
\par
Let $\Om \subseteq \R^3$ be an elastic body, $\partial_D \Om$ a part of its
boundary and let $g: \partial_D \Om \to \R^3$ be the spatial displacement of $\Om$ at
the points of $\partial_D \Om$. According to Griffith theory,
given a preexisting crack $\Gamma_1 \subseteq \Omb$,
the new crack $\Gamma$ and the displacement $u: \Om \setminus \Gamma \to \R^3$
associated to $g$ at the equilibrium minimizes the following
elastic energy
\begin{equation}
\label{griffithenergy}
\Es(v,g,\Gamma):= \int_\Om \mu |Ev|^2 +\lambda |tr Eu|^2 \,dx+ \hs^2(\Gamma),
\end{equation}
among all cracks $\Gamma$ with $\Gamma_1 \subseteq \Gamma$ and all displacements
$v: \Om \setminus \Gamma \to \R^3$ with $v=g$ on $\partial_D \Om \setminus \Gamma$. Here $Ev$
denotes the symmetric part of
the gradient of $v$, $tr$ denotes the trace of the matrix,
and $\hs^2$ denotes the two dimensional Hausdorff measure.
Griffith criterion thus involves a competition between the {\it bulk energy} given by
$\int_\Om \mu|Eu|^2 +\lambda|tr Eu|^2\,dx$ and the {\it surface energy} given by $\hs^2(\Gamma)$.
The boundary condition is required only on $\partial_D \Om \setminus \Gamma$
because the displacement in a fractured region is supposed to be not transmitted.
We indicate by $\Es(g,\Gamma)$ the minimum value of
(\ref{griffithenergy}) among all $v: \Om \setminus \Gamma \to \R^3$
with $v=g$ on $\partial_D \Om \setminus \Gamma$.
\par
Suppose that the boundary displacement $g$ varies with the time $t \in [0,1]$. The quasi-static
evolution $t \to \Gamma(t)$ proposed in \cite{FM} requires that:
\begin{itemize}
\item[(1)] $\Gamma(t)$ is increasing in time, i.e., $\Gamma(t_1) \subseteq \Gamma(t_2)$
for all $0 \le t_1 \le t_2 \le 1$;
\item[{}]
\item[(2)] $\Es(g(t),\Gamma(t)) \le \Es(g(t),\Gamma)$ for all cracks $\Gamma$ such that
$\cup_{s<t} \Gamma(s) \subseteq \Gamma$;
\item[{}]
\item[(3)] the elastic energy $\Es(g(t),\Gamma(t))$ is absolutely continuous in time.
\end{itemize}
Condition $(1)$ stands for the irreversibility of the evolution (fracture can only
increase); condition $(2)$ states that each time $t$ is of equilibrium, while condition
$(3)$ requires a regularity for the total energy.
\par
The problem of giving a precise mathematical formulation of the preceding model has been
the object of several recent papers. In 2000, Dal Maso and Toader \cite{DMT} dealt with
the case of antiplanar shear in dimension two: the authors consider a cylindric elastic body
$\Om= \Om' \times \R$ with $\Om' \subseteq \R^2$ subject to displacements of the form
$u(\pi x)e_3$ where $e_3$ is the versor of the $x_3$-axis, and $\pi$ is the projection on $\Om'$.
The boundary antiplanar displacement is assigned on $\partial_D \Om' \times \R$ while
the admissible cracks are of the form $K \times \R$ with $K$
compact connected subset of $\overline{\Om'}$ with a prescribed number of connected
components and with finite $\hs^1$-length. A generalization to non-isotropic surface energies
is contained in \cite{G}.
\par
Recently Francfort and Larsen \cite{FL} proposed a mathematical formulation which involves
the space $SBV$ of special functions of bounded variation. Their approach permits to treat
antiplanar shear in a $N$-dimensional setting, and allows to consider fractures with a possibly
infinite number of connected components. To be precise, they consider deformations of the
form $u(x) e_{N+1}$, where $u \in SBV(\Om)$ and $e_{N+1}$ denotes the unitary vector of the
$(N+1)$-axis.
The crack at time $t$ is defined as $\Gamma(t) \times \R$ where
$\Gamma(t):=\cup_{s<t} \left[ S_{u(s)} \cup
(\partial_D \Om \cap \{u(s) \not= g(s)\})
\right]$, and the pair $(u(t), \Gamma(t))$
is such that:
\begin{itemize}
\item[(a)]
for all $z \in SBV(\Om)$
\begin{equation}
\label{minintrprop}
\int_\Om |\nabla u(t)|^2 + \hs^{N-1}(\Gamma(t)) \le
\int_\Om |\nabla z|^2 + \hs^{N-1}((S_z \cup
(\partial_D \Om \cap \{z \not= g(t)\})
\cup \Gamma(t));
\end{equation}
\item[{}]
\item[(b)]
the elastic energy $\Es(t):= \int_\Om |\nabla u(t)|^2 + \hs^{N-1}(\Gamma(t))$ is absolutely
continuous and
$$
\Es(t)=\Es(0) + 2\int_0^t \int_\Om \nabla u(\tau) \nabla \dot{g}(\tau) \,dx \,d\tau.
$$
\end{itemize}
\indent
Numerical computations concerning this model of evolution (see \cite{BFM}) are
performed using a discretization in time procedure and an approximation of the elastic
energy proposed in 1990 by Ambrosio and Tortorelli (see \cite{AT1},\cite{AT2}).
Being the new energy elliptic, the difficulties arising in the discretization of the free
discontinuity term given by the fracture are avoided.
Supposing to have determined the deformation $u_i$ and the fracture $K_i$ at the time $t_i$, one
minimizes the Ambrosio-Tortorelli functional in the domain $\Om \setminus K_i$ under the
boundary conditions $g(t_{i+1})$, and hence reconstruct the couple $(u_{i+1},K_{i+1})$. In this way,
errors due to the discretization in time and to the approximation of the energy are introduced.
In order to study the convergence of the procedure, one is led to formulate a natural
notion of quasi-static evolution for the Ambrosio-Tortorelli functional.
The aim of this paper is to prove the convergence of this regular evolution
to an evolution of brittle fractures in the sense of \cite{FL}.
\par
The Ambrosio-Tortorelli functional is given by
\begin{equation*}
F_\varepsilon (u,v)= \ateps{u}{v}
\end{equation*}
where $(u,v) \in H^1(\Om) \times H^1(\Om)$, $0 \le v \le 1$,
$0 < \eta_\varepsilon << \varepsilon$.
$F_\varepsilon$ contains an {\it elliptic part}
\begin{equation}
\label{ellintr}
\elleps{u}{v}
\end{equation}
and a {\it surface part}
\begin{equation}
\label{surfintr}
MM_{\varepsilon}(v):=\mmeps{v}
\end{equation}
which is a term of Modica-Mortola type (see \cite{MoMo}).
\par
If a sequence $(\ueps,\veps)$ is such that $\Ateps{\ueps}{\veps} +||\ueps||_{\infty} \le C$,
then $\veps \to 1$
strongly in $L^2(\Om)$, and it turns out that, up to a subsequence, $\ueps \to u \in SBV(\Om)$;
roughly speaking, the gradient of $\ueps$ becomes larger and larger in the thick
regions in which $\veps$ approaches zero, possibly creating some jumps in the limit.
We conclude that the function $\ueps$ has to be considered as a regularization of the deformation
$u$, while the function $\veps$ has to be intended as a function which tends to $0$ in
the region where $S_u$ will appear, and to $1$ elsewhere.
Moreover (\ref{ellintr}) and (\ref{surfintr}) have to be
interpreted as regularizations of the bulk and surface elastic energy of $u$.
\par
In the regular context of the Ambrosio and Tortorelli functional, we define through a
variational argument the following notion of quasi-static evolution: we find a map
$t \to (u(t),v(t))$ from $[0,1]$ to $H^1(\Om) \times H^1(\Om)$, $0 \le v(t) \le 1$,
$u(t)=g(t)$, $v(t)=1$ on $\partial_D \Om$ such that:
\begin{itemize}
\item[(a)]
for all $0 \le s<t \le 1$: $v(t) \le v(s)$;
\item[{}]
\item[(b)]
for all $(u,v) \in H^1(\Om) \times H^1(\Om)$ with $u=g(t)$, $v=1$ on $\partial_D \Om$,
$0 \le v \le v(t)$:
\begin{equation}
\label{minatprop}
\Ateps{u(t)}{v(t)} \le \Ateps{u}{v};
\end{equation}
\item[{}]
\item[(c)]
the energy $\Es_\varepsilon(t):=\Ateps{u(t)}{v(t)}$ is absolutely continuous and for all
$t \in [0,1]$
$$
\Es_\varepsilon(t)=\Es_\varepsilon(0)+ 2\int_0^t \int_\Om (\eta_\varepsilon +v^2(\tau))
\nabla u(\tau) \nabla \dot{g}(\tau) \,dx \,d\tau;
$$
\item[{}]
\item[(d)]
there exists a constant $C$ depending only on $g$ such that $\Es_\varepsilon(t) \le C$ for all
$t \in [0,1]$.
\end{itemize}
Condition $(a)$ permits to recover in this regular context the fact that the fracture is
increasing in time: in fact, as $v(t)$ determines the fracture in the regions where it is
near zero, the condition $v(t) \le v(s)$ ensures that existing cracks are preserved at
subsequent times. Condition $(b)$ reproduces the minimality condition at each time with respect
to larger fractures, while condition $(c)$ describes the evolution in time of the total energy.
Condition $(d)$ gives the necessary compactness in order to let $\varepsilon \to 0$.
In the particular case in which $||g(t)||_{\infty} \le C_1$ for all $t \in [0,1]$, it
turns out that, using truncation arguments, $||\ueps(t)||_{\infty} \le C_1$ for all $t$ so that
a uniform $L^\infty$ bound is available at any time.
The requirement $v(t)=1$ on $\partial_D \Om$ for all $t \in [0,1]$ is
made in such a way that, letting $\varepsilon \to 0$, the surface energy
of the fracture in the limit is the usual one also for the part touching
the boundary $\partial_D \Om$.
\par
The main result of the paper is that, as $\varepsilon \to 0$, the quasi-static evolution
$t \to (\ueps(t),\veps(t))$ for the Ambrosio-Tortorelli functional converges to a quasi-static
evolution for brittle fracture in the sense of \cite{FL}. More precisely, there exists a
quasi-static evolution $t \to u(t) \in SBV(\Om)$ relative to the boundary data $g$ and a
sequence $\varepsilon_n \to 0$ such that for all $t \in [0,1]$ which are not discontinuity
points of $\hs^{N-1}(\Gamma(\cdot))$ we have
$$
v_n(t) \nabla \uepsn(t) \to \nabla u(t) \quad {\rm strongly\;in}\; L^2(\Om, \R^N),
$$
$$
\ellepsn{\uepsn(t)}{\vepsn(t)} \to \int_\Om |\nabla u(t)|^2 \, dx
$$
and
$$
MM_{\varepsilon_n}(\vepsn(t)) \to \hs^{N-1}(\Gamma(t)).
$$
Moreover $\Es_{\varepsilon_n}(t) \to \Es(t)$ for all $t \in [0,1]$.
We thus obtain an approximation of the total energy at any time, and an approximation of the
gradient of the deformation, of the bulk and the surface energy at all time
up to a countable set.
The main step in the proof is to derive the unilateral minimality
property~\eqref{minintrprop} from its regularized version~\eqref{minatprop}.
Given $z \in SBV(\Om)$, a natural way consists in constructing
$z_n \in H^1(\Om)$ and $v_n \in H^1(\Om)$ with $z_n=g(t)$, $v_n=1$ on $\partial_D \Om$,
$0 \le v_n \le v_n(t)$ and such that
\begin{equation}
\label{intrapprox1}
\lim_n \int_{\Om} (\eta_{\varepsilon_n}+v_n^2) |\nabla z_n|^2 \,dx=
\int_{\Om} |\nabla z|^2 \,dx
\end{equation}
and
\begin{equation}
\label{intrapprox2}
\limsup_n
\left[
MM_{\varepsilon_n}(v_n)-MM_{\varepsilon_n}(v_n(t))
\right]
\le \hs^{N-1} \left( S_z \setminus \Gamma(t) \right).
\end{equation}
We thus need a recovery sequence both for the deformation and the fracture: moreover
we have to take into account the boundary conditions and the constraint
$v_n \le v_n(t)$. Density results on $z$, such that of considering $S_z$ polyhedral,
cannot be directly applied since the set $S_z \cap \Gamma(t)$ could increase too
much; on the other hand it is not possible to work in $\Om \setminus \Gamma(t)$ since
no regularity results are available for $\Gamma(t)$ apart from its rectifiability.
It turns out that $S_z \cap \Gamma(t)$ is the part of the fracture more difficult
to be regularized, and in fact all the problems in the construction of $(z_n,v_n)$
are already present in the particular case $S_z \subseteq \Gamma(t)$.
In order to fix ideas, let us suppose
to be in this situation; we solve the problem in two steps. We firstly construct
$\tilde{z}_n \in SBV(\Om)$ with $\nabla \tilde{z}_n \to \nabla z$ strongly in
$L^2(\Om;\R^N)$ and such that $S_{\tilde{z}_n}$ is related to $u_n(t)$ and
$v_n(t)$ with precise energy estimates: this is done following the ideas of
\cite[Theorem 2.1]{FL}, that is using local reflections and gluing along the
boundaries of suitable upper levels of $u_n(t)$, but we have to choose the upper
levels in a more accurate way. In a second time, we regularize $S_{\tilde{z}_n}$
using not only $v_n(t)$, which is quite natural, but also $u_n(t)$, so that
\eqref{intrapprox1} and \eqref{intrapprox2} hold.
\par
The plan of the paper is the following. We introduce in Section \ref{notation}
the notation and the main tools employed in the rest of the paper.
In Section \ref{Atevol} we treat the quasi-static evolution
for the Ambrosio-Tortorelli functional, while in Section \ref{evofracture} we prove
the main approximation result. The derivation of the minimality
property~\eqref{minintrprop} is contained in Section \ref{secminthm}.

\section{Notation and Preliminaries}
\label{notation}

In this section we state the notations and introduce the main tools
used in the rest of the paper.

\vskip10pt
{\it Basic notation.}
In the rest of the paper, we will employ the following basic notations:
\begin{itemize}
\item[-] $\Om$ is an open bounded subset of $\R^N$ with Lipschitz boundary;
\item[-] $\partial_D \Om$ is a subset on $\partial \Om$
open in the relative topology;
\item[-] $L^p(\Om;\R^m)$ with $1 \le p \le +\infty$ and $m \ge 1$ is the
Lebesgue space of $p$-summable $\R^m$-valued functions;
\item[-] $H^1(\Om)$ is the Sobolev spaces of functions in $L^2(\Om)$
with distributional derivative in $L^2(\Om;\R^N)$;
\item[-] if $u \in H^1(\Om)$, $\nabla u$ is its gradient;
\item[-] $\hs^{N-1}$ is the $(N-1)$-dimensional Hausdorff measure;
\item[-] $\|\cdot\|_{\infty}$ is the sup-norm;
\item[-] $1_A$ is the characteristic function of $A$;
\item[-] if $\sigma \in ]0,+\infty[$, $o(\sigma)$ is such that
$\lim_{\sigma \to 0^+} o(\sigma)=0$.
\end{itemize}

\vskip10pt
{\it Special functions of bounded variation.}
For the general theory of functions of bounded variation, we refer to \cite{AFP}; here
we recall some basic definitions and theorems we need in the sequel.
\par
Let $A$ be an open subset of $\R^N$, and let $u: A \to \R^n$.
We say that $u \in BV(A;\R^n)$ if $u \in L^1(A;\R^n)$, and its distributional
derivative is a vector-valued Radon measure on $A$.
\par
We say that $u \in SBV(A;\R^n)$ if $u \in BV(A;\R^n)$ and its distributional
derivative can be represented as
$$
Du(B)= \int_B \nabla u(x) \,dx+ \int_{B \cap S_u} (u^+(x)-u^-(x)) \otimes \nu_x
\,d\hs^{N-1}(x)
$$
where $\nabla u$ denotes the approximate gradient of $u$, $S_u$ denotes the set of
approximate jumps of $u$, $u^+$ and $u^-$ are the traces of $u$ on $S_u$, and $\nu_x$
is the normal to $S_u$ at $x$. The space $SBV(A;\R^n)$ is called the space of
special functions of bounded variation. Note that if $u \in SBV(A;\R^n)$, then
the singular part of $Du$ is concentrated on $S_u$ which turns out to be countably
$\hs^{N-1}$-rectifiable.
\par
We say that $u \in GSBV(A)$ if for every $M \ge 0$ we have
$(u \wedge M) \vee (-M) \in SBV(A)$. For every $p \in [1, \infty]$, let us pose
$SBV^p(A, \R^n):=\{u \in SBV(A, \R^n)\,:\, \nabla u \in L^p(A; M^{N \times n}) \}$, and
$GSBV^p(A):= \{u \in GSBV(A)\,:\, \nabla u \in L^p(A; \R^N) \}$.
\par
The space $SBV$ is very useful when dealing with variational problems involving
volume and surface energies because of the following compactness and lower
semicontinuity result due to L.Ambrosio (see \cite{A1}, \cite{A3}).

\begin{theorem}
\label{SBVcompact}
Let $(u_k)$ be a sequence in $SBV(A;\R^n)$ such that there exist $q>1$ and $c \ge 0$
with
$$
\int_A |\nabla u_k|^q \,dx+ \hs^{N-1}(S_{u_k})+ ||u_k||_{\infty,A} \le c
$$
for every $k \in \N$. Then there exist a subsequence $(u_{k_h})$ and a function
$u \in SBV(A;\R^n)$ such that
\begin{align}
\label{sbvconv}
\nonumber
u_{k_h} \to u \quad {strongly \; in}\; L^1(A;\R^n), \\
\nabla u_{k_h} \weak \nabla u \quad {weakly \; in}\; L^1(A;M^{N \times n}), \\
\nonumber
\hs^{N-1}(S_u) \le \liminf_h \hs^{N-1}(S_{u_{k_h}}).
\end{align}
\end{theorem}

In the rest of the paper, we will say that $u_k \to u$ in $SBV(\Om)$ if
$u_k$ and $u$ satisfy \eqref{sbvconv}.

\vskip10pt
{\it Quasi-static evolution of brittle fractures.}
Let $g:[0,1] \to H^1(\Om)$ be absolutely continuous; we indicate
the gradient of $g$ at time $t$ by $\nabla g(t)$, and the time
derivative of $g$ at time $t$ by $\dot{g}(t)$.
The main result of \cite{FL} is the following theorem.

\begin{theorem}
\label{qse}
There exists a crack $\Gamma(t) \subseteq \Omb$ and a field $u(t) \in SBV(\Om)$ such
that
\begin{itemize}
\item[(a)] $\Gamma(t)$ increases with $t$;
\item[{}]
\item[(b)] $u(0)$ minimizes
$$
\int_\Om |\nabla v|^2 \,dx +
\hs^{N-1}( S_v \cup
\{x \in \partial_D \Om: v(x) \not= g(0)(x)\})
$$
among all $v \in SBV(\Om)$ (inequalities on $\partial_D \Om$ are intended for the
traces of $v$ and $g$);
\item[{}]
\item[(c)]
for $t>0$, $u(t)$ minimizes
$$
\int_\Om |\nabla v|^2 \,dx
+\hs^{N-1} \left( \left[ S_v \cup
\{x \in \partial_D \Om: v(x) \not= g(t)(x)\} \right] \setminus \Gamma(t)
\right)
$$
among all $v \in SBV(\Om)$;
\item[{}]
\item[(d)]
$S_{u(t)} \cup \{x \in \partial_D \Om\,:\, u(t)(x) \not= g(t)(x) \}
\subseteq \Gamma(t)$, up to a set of $\hs^{N-1}$--measure $0$.
\end{itemize}
Furthermore, the total energy
$$
\Es(t):= \int_\Om |\nabla u(t)|^2 \,dx +\hs^{N-1}( \Gamma(t))
$$
is absolutely continuous and is given by
$$
\Es(t)=\Es(0)+2 \int_0^t \int_\Om \nabla u(\tau) \nabla \dot{g}(\tau) \,dx \,d\tau.
$$
Finally, for any countable, dense set $I \subseteq [0,1]$, the crack $\Gamma(t)$
and the field $u(t)$ can be chosen such that
$$
\Gamma(t)= \bigcup_{\tau \in I, \tau \le t}
\left( S_{u(\tau)} \cup
\{x \in \partial_D \Om\,:\, u(\tau)(x) \not= g(\tau)(x) \} \right)
$$
\end{theorem}

\vskip10pt
{\it The Ambrosio-Tortorelli functional.}
In \cite{AT1} and \cite{AT2}, Ambrosio and Tortorelli proposed an elliptic
approximation of the Mumford-Shah functional in the sense of $\Gamma${-}convergence.
Their result has been extended in the vectorial case
in \cite{Fo}, where non-isotropic surface energies are also considered.
For every $u \in GSBV(\Om)$ let
$$
F(u):=\int_\Om |\nabla u|^2 \,dx +\hs^{N-1}(S_u)
$$
the well known Mumford-Shah functional; for every
$(u,v) \in H^1(\Om) \times H^1(\Om)$ the Ambrosio-Tortorelli functional is defined by
\begin{equation*}
\Ateps{u}{v}:= \ateps{u}{v}
\end{equation*}
where $\eta_\varepsilon >0$ and $\eta_\varepsilon << \varepsilon$.
Let us indicate the space of Borel functions on $\Om$ by $\bs(\Om)$ and let us
consider on $\bs(\Om) \times \bs(\Om)$ the functionals
$$
\fs(u,v,\Om):=
\left\{
\begin{array}{l}
F(u) \\
{}\\
+\infty
\end{array}
\begin{array}{l}
u \in GSBV(\Om), v \equiv 1 \;{\rm a.e. \;on}\; \Om \\
{}\\
{\rm otherwise}
\end{array}
\right.
$$
and
$$
\fs_\varepsilon(u,v,\Om):=
\left\{
\begin{array}{l}
\Ateps{u}{v} \\
{}\\
+\infty
\end{array}
\begin{array}{l}
(u,v) \in H^1(\Om), 0 \le v \le 1 \\
{}\\
{\rm otherwise.}
\end{array}
\right.
$$

The Ambrosio-Tortorelli result can be expressed in the following way.

\begin{theorem}
\label{at}
The functionals $(\fs_\varepsilon)$ on $\bs(\Om) \times \bs(\Om)$
$\Gamma${-}converge to $\fs$ with respect to the convergence in measure.
\end{theorem}

In particular, we will use several times the following fact: if $(\ueps,\veps) \in H^1(\Om)
\times H^1(\Om)$ is such that $F_{\varepsilon}(\ueps,\veps) +||\ueps||_{\infty} \le C$, there
exists $u \in SBV(\Om)$ and a sequence $\varepsilon_k \to 0$ such that $u_{\varepsilon_k} \to
u$ a.e., and
\begin{equation}
\label{atliminf}
\int_\Om |\nabla u|^2 \,dx +\hs^{N-1}(S_u) \le \liminf_{\varepsilon \to 0}
F_{\varepsilon}(\ueps,\veps).
\end{equation}

\vskip10pt
{\it A density result.}
Let $A \subseteq \R^N$ be open.
We say that $K \subseteq A$ is polyhedral (with respect to $A$),
if it is the intersection of $A$ with the union of a finite number of $(N-1)$-dimensional
simplexes of $S$.
\par
The following density result is proved in \cite{C}.

\begin{theorem}
\label{piecedensity}
Assume that $\partial A$ is locally Lipschitz, and let $u \in GSBV^p(A)$. For
every $\varepsilon>0$, there exists a function $v \in SBV^p(A)$ such that
\begin{itemize}
\item[(a)]
$S_v$ is essentially closed, i.e., $\hs^{N-1}( \overline{S_v} \setminus S_v)=0$;
\item[{}]
\item[(b)]
$\overline{S_v}$ is a polyhedral set;
\item[{}]
\item[(c)]
$v \in W^{k, \infty}(A \setminus \overline{S_v})$ for every $k \in \N$;
\item[{}]
\item[(d)]
$||v-u||_{L^p(A)} < \varepsilon$;
\item[{}]
\item[(e)]
$||\nabla v- \nabla u||_{L^p(A; \R^N)} < \varepsilon$;
\item[{}]
\item[(f)]
$|\hs^{N-1}(S_v)-\hs^{N-1}(S_u)| <\varepsilon$.
\end{itemize}
\end{theorem}

Theorem \ref{piecedensity} has been generalized to non-isotropic surface energies in
\cite{CT}. In Section \ref{secminthm}, we will use the following result.

\begin{proposition}
\label{regularization}
Let $\partial_N \Om:= \partial \Om \setminus \partial_D \Om$,
$B$ an open ball such that $\Omb \subseteq B$, and let
$\Om':=B \setminus \partial_N \Om$. Given $g \in H^1(B)$ and $u \in SBV(\Om')$
with $u=g$ on $\Om' \setminus \Omb$, there exists $u_h \in SBV(\Om')$ such that
\begin{itemize}
\item[(a)] $u_h=g$ in $\Om' \setminus \Omb$ and in a neighborhood of
$\partial_D \Om$;
\item[{}]
\item[(b)] $\overline{S_{u_h}}$ is polyhedral and $\overline{S_{u_h}} \subseteq \Om$ for all $h$;
\item[{}]
\item[(c)] $\nabla u_h \to \nabla u$ strongly in $L^2(\Om';\R^N)$;
\item[{}]
\item[(d)] for all $A$ open subset of $\Om'$ with $\hs^{N-1}(\partial A \cap S_u)=0$,
we have
$$
\lim_h \hs^{N-1}(A \cap S_{u_h})=\hs^{N-1}(A \cap S_u).
$$
\end{itemize}
\end{proposition}

\begin{proof}
Using a partition of unity, we may prove the result in the case
$\Om':=Q \times ]-1,1[$, $\Om:=\{(x,y) \in Q \times ]-1,1[\,:\, y >f(x)\}$,
$\partial_D \Om:=\{(x,y) \in Q \times ]-1,1[\,:\,y=f(x)\}$, where $Q$ is unit cube
in $\R^{N-1}$ and $f:Q \to \R$ is a Lipshitz function with values in $]-\frac{1}{2},\frac{1}{2}[$.
Let $g \in H^1(\Om')$, and let $u \in SBV(\Om')$ with $u=g$ on
$\Om' \setminus \Om$.
\par
Let $w_h:=u(x-he_N)$ where $e_N$ is the versor of the $N$-axis, and let $\varphi_h$
be a cut off function with $\varphi_h=1$ on $\{y \le f(x)+\frac{h}{3}\}$,
$\varphi_h=0$ on $\{y \ge f(x)+\frac{h}{2}\}$, and $||\nabla \varphi_h||_\infty
\le \frac{1}{h}$. Let us pose $v_h:= \varphi_h g+ (1-\varphi_h)w_h$.
We have that $v_h=g$ in $\Om' \setminus \Omb$ and in a neighborhood of
$\partial_D \Om$; moreover we have
$$
\nabla v_h= \varphi_h \nabla g+ (1-\varphi_h) \nabla w_h+\nabla \varphi_h(g-w_h).
$$
Since $\nabla \varphi_h (g-w_h) \to 0$ strongly in $L^2(\Om';\R^N)$, we have
$\nabla v_h \to \nabla u$ strongly in $L^2(\Om';\R^N)$. Finally, for all
$A$ open subset of $\Om'$ with $\hs^{N-1}(\partial A \cap S_u)=0$,
we have
$$
\lim_h \hs^{N-1}(A \cap S_{v_h})=\hs^{N-1}(A \cap S_u).
$$
\par
In order to conclude the proof, let us apply Theorem \ref{piecedensity}
obtaining $\tilde{v}_h$ with polyhedral jumps such that
$||v_h-\tilde{v}_h||_{L^2(\Om')}+
||\nabla v_h -\nabla \tilde{v}_h||_{L^2(\Om';\R^N)} \le h^2$,
$|\hs^{N-1}(S_{v_h}) -\hs^{N-1}(S_{\tilde{v}_h})| \le h$. If we pose
$u_h:= \varphi_h g+(1-\varphi_h) \tilde{v}_h$, we obtain the thesis.
\end{proof}

\section{The Main Results}
\label{mainresult}

If $g \in W^{1,1}([0,1];H^1(\Om))$, we indicate
the gradient of $g$ at time $t$ by $\nabla g(t)$, and the time
derivative of $g$ at time $t$ by $\dot{g}(t)$.
\par
Concerning the Ambrosio-Tortorelli functional, the following theorem holds.

\begin{theorem}
\label{atevol}
Let $g \in W^{1,1}([0,1];H^1(\Om))$. Then there exists a strongly measurable map
$$
\begin{array}{c}
[0,1] \\ t
\end{array}
\begin{array}{c}
\longrightarrow \\ \longmapsto
\end{array}
\begin{array}{c}
H^1(\Om) \times H^1(\Om) \\ (u(t),v(t))
\end{array}
$$
such that $0 \le v(t) \le 1$ in $\Om$, $u(t)=g(t)$, $v(t)=1$ on $\partial_D \Om$
for all $t \in [0,1]$, and:
\begin{itemize}
\item[(a)] for all $0 \le s \le t \le 1: v(t) \le v(s)$;
\item[{}]
\item[(b)] for all $(u,v) \in H^1(\Om) \times H^1(\Om)$ with
$u=g(0)$, $v=1$ on $\partial_D \Om$
$$
\Ateps{u(0)}{v(0)} \le \Ateps{u}{v};
$$
\item[{}]
\item[(c)] for all $t \in ]0,1]$ and for all $(u,v) \in H^1(\Om) \times H^1(\Om)$
with $0 \le v \le v(t)$ on $\Om$, and $u=g(t)$, $v=1$ on $\partial_D \Om$
$$
\Ateps{u(t)}{v(t)} \le \Ateps{u}{v};
$$
\item[{}]
\item[(d)] the function $t \to \Ateps{u(t)}{v(t)}$ is absolutely continuous
and
$$
\Ateps{u(t)}{v(t)}=\Ateps{u(0)}{v(0)}+ 2
\int_0^t \int_{\Om} (\eta_\varepsilon+v^2(\tau))
\nabla u(\tau) \nabla \dot{g}(\tau) \,dx\,d\tau.
$$
\end{itemize}
\end{theorem}

The main result of the paper is the following theorem.

\begin{theorem}
\label{mainthm}
Let $g \in W^{1,1}([0,1];H^1(\Om))$ such that $\|g(t)\|_\infty \le C$ for all
$t \in [0,1]$, and let  $g_h \in W^{1,1}([0,1];H^1(\Om))$ be a sequence of absolutely
continuous functions with $\|g_h(t)\|_\infty \le C$, $g_h(t) \in C(\Omb)$ for all
$t \in [0,1]$ and such that
$g_h \to g$ strongly in $W^{1,1}([0,1];H^1(\Om))$.
For all $\varepsilon>0$, let $t \to (u_{\varepsilon,h}(t),v_{\varepsilon,h}(t))$
be a quasi-static
evolution for the Ambrosio-Tortorelli functional $F_\varepsilon$ with
boundary data $g_h$ given by Theorem \ref{atevol}.
\par
Then there exists a quasi-static evolution $t \to u(t) \in SBV(\Om)$
relative to the boundary data $g$ in the sense of Theorem \ref{qse},
and two sequences $\varepsilon_n \to 0$ and $h_n \to +\infty$ such that, posing
$u_n:= u_{\varepsilon_n,h_n}$ and $v_n:=v_{\varepsilon_n,h_n}$, the following hold:
\begin{itemize}
\item[(a)]
for all $t \in [0,1]$ we have
\begin{equation*}
F_{\varepsilon_n}(u_n(t),v_n(t)) \to \Es(t);
\end{equation*}
\item[{}]
\item[(b)] if $\ns$ denotes the point of discontinuity of
$\hs^{N-1}(\Gamma(\cdot))$, for all $t \in [0,1] \setminus \ns$ we have
\begin{equation*}
v_n(t) \nabla u_n(t) \to \nabla u(t) \quad {\rm strongly \; in}\;L^2(\Om;\R^N),
\end{equation*}
\begin{equation*}
\lim_n \int_{\Om} (\eta_n+v^2_n(t)) |\nabla u_n(t)|^2 \,dx = \int_{\Om}
|\nabla u(t)|^2 \,dx,
\end{equation*}
and
\begin{equation*}
\lim_n \mmn{v_n(t)} =\hs^{N-1}(\Gamma(t)).
\end{equation*}
\end{itemize}
\end{theorem}

Theorem \ref{atevol} concerning the quasi-static evolution for the Ambrosio-Tortorelli
functional is proved in Section \ref{Atevol}.
In Section \ref{evofracture} we prove the compactness and approximation result
given by Theorem \ref{mainthm}. An important step in the proof is given by Theorem
\ref{minthm} to which is dedicated the entire Section \ref{secminthm}.

\section{Quasi-static evolution for the Ambrosio-Tortorelli functional}
\label{Atevol}

This section is devoted to the proof of Theorem \ref{atevol} where
a suitable notion of quasi-static evolution in a regular context is
proposed. The evolution will be obtained through a discretization in time
procedure: each step will be performed using a variational argument
which will give the minimality property stated in points $(b)$ and $(c)$.
\par
Let $g \in W^{1,1}([0,1];H^1(\Om))$.
Given $\delta>0$, let $N_\delta$ be the largest integer such that
$\delta N_\delta \le 1$; for $i \ge 0$ we pose $t_i^\delta=i\delta$ and
for $0 \le i \le N_\delta$ we pose $g_i^\delta=g(t_i^\delta)$. Define
$u^\delta_0$ and $v^\delta_0$ as a minimum for the problem
\begin{equation}
\label{atstep1}
\min \{\Ateps{u}{v}\,:\,(u,v) \in H^1(\Om) \times H^1(\Om),
0 \le v \le 1 \;{\rm in}\; \Om,
u=g_0^\delta, v=1 \,{\rm on}\, \partial_D \Om\},
\end{equation}
and let $(u_{i+1}^\delta, v_{i+1}^\delta)$ be a minimum
for the problem
\begin{equation}
\label{atnextstep}
\min \{\Ateps{u}{v}\,:\,(u,v) \in H^1(\Om) \times H^1(\Om),
0\le v \le v_i^\delta \;{\rm in}\; \Om,
u=g_{i+1}^\delta, v=1 \,{\rm on}\, \partial_D \Om\}.
\end{equation}
Problems \eqref{atstep1} and \eqref{atnextstep} are well posed: in fact,
referring for example to problem \eqref{atnextstep}, let $(u_n,v_n)$ be
a minimizing sequence. Since $(g_{i+1}^\delta, v_{i}^\delta)$ is an
admissible pair, we obtain that there exists a constant $C>0$ such that
for all $n$
$$
\Ateps{u_n}{v_n} \le C.
$$
Since $\varepsilon, \eta_\varepsilon>0$, we deduce that $(u_n,v_n)$
is bounded in $H^1(\Om) \times H^1(\Om)$ so that up to a subsequence
$u_n \weak u$ and $v_n \weak v$ weakly in $H^1(\Om)$. We get immediately
that $u=g_{i+1}^\delta$ and $v=1$ on $\partial_D \Om$ since
$u_n=g_{i+1}^\delta$ and $v_n=1$ on $\partial_D \Om$
for all $n$; on the other hand, since $v_n \to v$ strongly in $L^2(\Om)$,
we obtain that $0 \le v \le v_i^\delta$. By semicontinuity, we have
$$
\Ateps{u}{v} \le \liminf_n \Ateps{u_n}{v_n}
$$
so that $(u,v)$ is a minimum point for problem (\ref{atnextstep}).
\par
We note that by minimality of the pair $(u_{i+1}^\delta,v_{i+1}^\delta)$, we may
write
\begin{multline}
\label{energyestimate}
\Ateps{u_{i+1}^\delta}{v_{i+1}^\delta} \le
\Ateps{u_i^\delta+g_{i+1}^\delta-g_i^\delta}{v_i^\delta}= \\
=\Ateps{u_i^\delta}{v_i^\delta}
+2\int_\Om (\eta_\varepsilon + (v_i^\delta)^2)
\nabla u_i^\delta \nabla(g_{i+1}^\delta -g_i^\delta) \,dx
+\int_\Om (\eta_\varepsilon + (v_i^\delta)^2)
|\nabla(g_{i+1}^\delta -g_i^\delta)|^2 \,dx \le \\
\le
\Ateps{u_i^\delta}{v_i^\delta}
+2\int_{t_i^\delta}^{t_{i+1}^\delta} \int_\Om (\eta_\varepsilon + (v_i^\delta)^2)
\nabla u_i^\delta \nabla \dot{g}(\tau) \,dx \,d\tau
+e(\delta) \int_{t_i^\delta}^{t_{i+1}^\delta}
||\nabla \dot{g}(\tau)||_{L^2(\Om;\R^N)}\,d\tau,
\end{multline}
where
$$
e(\delta):= (1+\eta_{\varepsilon})
\max_{0 \le r \le N_\delta -1} \int_{t_r^\delta}^{t_{r+1}^\delta}
||\nabla \dot{g}(\tau)||_{L^2(\Om;\R^N)} \,d\tau
$$
is infinitesimal as $\delta \to 0$.
\par
We now make a piecewise constant interpolation defining
\begin{equation}
\label{discrevol}
u^\delta(t)=u_i^\delta, \quad v^\delta(t)=v_i^\delta, \quad
g^\delta(t)=g_i^\delta \quad \mbox{ for }t_i^\delta \le t <t_{i+1}^\delta.
\end{equation}
Note that
by construction the map $t \to v^\delta(t)$ is decreasing from $[0,1]$ to $L^2(\Om)$.
Moreover, iterating the estimate (\ref{energyestimate}), we obtain
\begin{eqnarray}
\label{mainestimate}
\nonumber
\Ateps{u^\delta(t)}{v^\delta(t)} &\le&
\Ateps{u^\delta(s)}{v^\delta(s)}
+2\int_{s^\delta}^{t^\delta} \int_{\Om} (\eta_\varepsilon + v^\delta(\tau)^2)
\nabla u^\delta(\tau) \nabla \dot{g}(\tau) \,dx \,d\tau +\\
&&+e(\delta) \int_{s^\delta}^{t^\delta}||\nabla \dot{g}(\tau)||_{L^2(\Om;\R^N)}\,d\tau
\end{eqnarray}
where $s^\delta:=s_i^\delta$ with $s_i^\delta \le s <s_{i+1}^\delta$,
$t^\delta:=t_i^\delta$ with $t_i^\delta \le t <t_{i+1}^\delta$.
\par
Note that by minimality of the pair $(u^\delta(t),v^\delta(t))$, we have
$$
\Ateps{u^\delta(t)}{v^\delta(t)} \le \Ateps{g^\delta(t)}{v^\delta(t)}
$$
so that
\begin{equation}
\label{estu}
\int_{\Om} (\eta_\varepsilon + v^\delta(t)^2) |\nabla u^\delta(t)|^2 \,dx
\le
\int_{\Om} (\eta_\varepsilon + v^\delta(t)^2) |\nabla g^\delta(t)|^2 \,dx \le C_1
\end{equation}
with $C_1>0$ independent of $\delta$ and $t$. In particular by (\ref{estu}) we have that
\begin{equation*}
||\nabla u^\delta(t)||^2_{L^2(\Om;\R^N)} \le \frac{C_1}{\eta_\varepsilon}.
\end{equation*}
Since $u^\delta(t)=g^\delta(t)$ on $\partial_D \Om$, and $g^\delta(t)$ is uniformly bounded
in $H^1(\Om)$ for all $t$ and $\delta$, we get by a variant of Poincar\'e inequality
that $u^\delta(t)$ is uniformly bounded in $H^1(\Om)$ for all $t$ and $\delta$.
\par
Now we pass to $v$ in order to obtain some coerciveness in the space $H^1(\Om)$.
Notice that
\begin{eqnarray*}
&2&\left| \int_0^{t^\delta} \int_{\Om} (\eta_\varepsilon + v^\delta(\tau)^2)
\nabla u^\delta(\tau) \nabla \dot{g}(\tau) \,dx \,d\tau \right| \le \\
\nonumber
&\le&
2\int_0^{t^\delta}\sqrt{\eta_\varepsilon+1}
\left(
\int_{\Om}(\eta_\varepsilon +v^\delta(t)^2)|\nabla u^\delta(t)|^2 \,dx
\right)^\frac{1}{2}
||\nabla \dot{g}(\tau)||_{L^2(\Om;\R^N)} \,d\tau,
\end{eqnarray*}
and by (\ref{estu}), we obtain
\begin{equation}
\label{estder}
\left| 2\int_0^{t^\delta} \int_{\Om} (\eta_\varepsilon + v^\delta(\tau)^2)
\nabla u^\delta(\tau) \nabla \dot{g}(\tau) \,dx \,d\tau \right| \le C_2
\end{equation}
with $C_2>0$ independent of $t$ and $\delta$.
\par
By (\ref{mainestimate}) with $s=0$, and (\ref{estder}), we deduce
\begin{eqnarray*}
&\displaystyle \frac{\varepsilon}{2}&\int_{\Om} |\nabla v^\delta(t)|^2 \,dx +
\frac{1}{2\varepsilon}\int_{\Om}(1-v^\delta(t))^2 \,dx \le \\
\nonumber
&\le& \Ateps{u^\delta(0)}{v^\delta(0)}
+2\int_{0}^{t^\delta} \int_{\Om} (\eta_\varepsilon + v^\delta(\tau)^2)
\nabla u^\delta(\tau) \nabla \dot{g}(\tau) \,dx \,d\tau + \\
\nonumber
&&
+e(\delta) \int_{s^\delta}^{t^\delta}||\nabla \dot{g}(\tau)||_{L^2(\Om;\R^N)}\,d\tau \le \\
\nonumber
&\le&
\Ateps{u^\delta(0)}{v^\delta(0)}+C_2
+e(\delta) \int_0^1||\nabla \dot{g}(\tau)||_{L^2(\Om} \,d\tau.
\end{eqnarray*}
We conclude that there exists $C>0$ independent of $t$ and $\delta$ such that for
all $t \in [0,1]$
\begin{equation}
\label{estvh1}
||v^\delta(t)||_{H^1(\Om)} \le C.
\end{equation}
We now want to pass to the limit in $\delta$ as $\delta \to 0$.

\begin{lemma}
\label{convv}
There exists a sequence $\delta_n \to 0$
and a strongly measurable map $v:[0,1] \to H^1(\Om)$ such that
$v^{\delta_n}(t) \weak v(t)$ weakly in $H^1(\Om)$ for all $t \in [0,1]$.
Moreover, $v$ is decreasing from $[0,1]$ to $L^2(\Om)$, and
$0 \le v(t) \le 1$ in $\Om$, $v(t)=1$ on $\partial_D \Om$ for all
$t \in [0,1]$.
\end{lemma}

\begin{proof}
Since the map $t \to v^\delta(t)$ is monotone decreasing from $[0,1]$ to $L^2(\Om)$,
and $0 \le v^\delta(t) \le 1$ for all $t$, we deduce
by a variant of Helly's compactness theorem for sequences of monotone real
functions, that there exists a subsequence $\delta_n \to 0$ and a decreasing
map $v:[0,1] \to L^2(\Om)$ such that for all $t \in [0,1]$ we have
$v^{\delta_n}(t) \to v(t)$ strongly in $L^2(\Om)$. In particular we deduce
$0 \le v(t) \le 1$ in $\Om$.
By (\ref{estvh1}), we have that for all $t \in [0,1]$, up to a subsequence, $v_{\delta_n}(t)
\weak w$ weakly in $H^1(\Om)$; since $v_{\delta_n}(t) \to v$ strongly in $L^2(\Om)$, we
deduce that $w=v(t)$ so that $v(t) \in H^1(\Om)$, and $v_{\delta_n}(t) \weak v(t)$ weakly in
$H^1(\Om)$. As a consequence, $v(t)=1$ on $\partial_D \Om$ for all $t \in [0,1]$.
Finally, $v$ is strongly measurable from $[0,1]$ to $H^1(\Om)$ because it is weakly
measurable and separably valued (see \cite[Chapter V, Section 4]{Yo}).
\end{proof}

Let us consider the sequence $\delta_n$, and the map $v$ given by Lemma \ref{convv}.
We indicate $u^{\delta_n}$, $v^{\delta_n}$ and $g^{\delta_n}$ simply by $u_n$, $v_n$
and $g_n$.

\begin{lemma}
\label{convu}
There exists a strongly measurable map $u:[0,1] \to H^1(\Om)$ such that
$u_n(t) \to u(t)$ strongly in $H^1(\Om)$ for all $t \in [0,1]$. In particular,
$u(t)=g(t)$ on $\partial_D \Om$ for all $t \in [0,1]$.
\end{lemma}

\begin{proof}
Let $t \in [0,1]$. We note that $u_n(t)$ is the minimum of the following problem
$$
\min \left\{ \int_{\Om} (\eta_\varepsilon +v_n^2(t)) |\nabla z|^2 \,dx
\,:\, z \in H^1(\Om),
z=g_n(t) \,{\rm on}\, \partial_D \Om \right\}.
$$
Since by Lemma \ref{convv} $v_n(t) \to v(t)$ strongly in $L^2(\Om)$, and
$g_n(t) \to g(t)$ strongly in $H^1(\Om)$, we deduce by standard results on
$\Gamma${-}convergence (see \cite{DM}),
that $u_n(t) \weak u(t)$ weakly in $H^1(\Om)$ where $u(t)$ is the solution of the problem
$$
\min \left\{ \int_{\Om} (\eta_\varepsilon +v^2(t)) |\nabla z|^2 \,dx
\,:\, z \in H^1(\Om),
z=g(t) \,{\rm on}\, \partial_D \Om \right\}.
$$
Moreover, we have also convergence of energies, that is
\begin{equation}
\label{minconv}
\lim_n \int_{\Om} (\eta_\varepsilon +v_n^2(t)) |\nabla u_n(t)|^2 \,dx=
\int_{\Om} (\eta_\varepsilon +v^2(t)) |\nabla u(t)|^2 \,dx.
\end{equation}
Since $v_n(t) \nabla u_n(t) \weak v(t) \nabla u(t)$ weakly in $L^2(\Om;\R^N)$,
we obtain
$$
\int_{\Om} v^2(t) |\nabla u(t)|^2 \,dx \le \liminf_n
\int_{\Om} v_n^2(t) |\nabla u_n(t)|^2 \,dx,
$$
so that by (\ref{minconv}) we deduce $\nabla u_n(t) \to \nabla u(t)$ strongly in
$L^2(\Om;\R^N)$. We conclude that $u_n(t) \to u(t)$ strongly in $H^1(\Om)$
for all $t \in [0,1]$, and so the map $t \to u(t)$ is strongly measurable from
$[0,1]$ to $H^1(\Om)$. Finally $u(t)=g(t)$ on $\partial_D \Om$ and the proof is
complete.
\end{proof}

The following minimality property for the pair $(u(t),v(t))$ holds.

\begin{proposition}
\label{minimality}
If $t \in ]0,1]$, for every $(u,v) \in H^1(\Om) \times H^1(\Om)$ such that
$0 \le v \le v(t)$ in $\Om$, and $u=g(t)$, $v=1$ on $\partial_D \Om$, we have
\begin{equation*}
\Ateps{u(t)}{v(t)} \le \Ateps{u}{v}.
\end{equation*}
Moreover, for all $(u,v) \in H^1(\Om) \times H^1(\Om)$ such that
$u=g(0)$, $v=1$ on $\partial_D \Om$, we have
\begin{equation*}
\Ateps{u(0)}{v(0)} \le \Ateps{u}{v}.
\end{equation*}
\end{proposition}

\begin{proof}
Let us pose
$$
u_n:= u+g_n(t)-g(t),
$$
and
$$
v_n:= \min \{v_n(t), v\};
$$
we have $u_n \to u$ strongly in $H^1(\Om)$, and $v_n \weak v$ weakly in $H^1(\Om)$.
Since $0 \le v_n \le v_n(t)$ in $\Om$, and $u_n=g_n(t)$, $v_n=1$ on $\partial_D \Om$,
by the minimality property of the pair $(u_n(t), v_n(t))$ we get
$$
\Ateps{u_n(t)}{v_n(t)} \le \Ateps{u_n}{v_n},
$$
that is
\begin{eqnarray}
\label{comp1}
\ateps{u_n(t)}{v_n(t)} \le \\
\nonumber
\le \ateps{u_n}{v_n}.
\end{eqnarray}
Notice that
$$
\frac{\varepsilon}{2} \int_{\Om} |\nabla v_n|^2 \,dx
=
\frac{\varepsilon}{2} \int_{\{v_n(t) < v\}} |\nabla v_n(t)|^2 \,dx
+\frac{\varepsilon}{2} \int_{\{v_n(t) \ge v\}} |\nabla v|^2 \,dx
$$
so that (\ref{comp1}) becomes
\begin{eqnarray*}
\nonumber
\elleps{u_n(t)}{v_n(t)}+
\frac{\varepsilon}{2} \int_{\{v_n(t) \ge v\}} |\nabla v_n(t)|^2 \,dx
+\frac{1}{2\varepsilon} \int_{\Om} (1-v_n(t))^2 \,dx \le \\
\le \elleps{u_n}{v_n}+
\frac{\varepsilon}{2} \int_{\{v_n(t) \ge v\}} |\nabla v|^2 \,dx
+\frac{1}{2\varepsilon} \int_{\Om} (1-v_n)^2 \,dx.
\end{eqnarray*}
For $n \to \infty$, the right hand side is less than $\Ateps{u}{v}$.
Let us consider the left hand side. By semicontinuity we have
\begin{equation*}
\liminf_n
\frac{\varepsilon}{2} \int_{\{v_n(t) \ge v\}} |\nabla v_n(t)|^2 \,dx
\ge
\frac{\varepsilon}{2} \int_\Om |\nabla v(t)|^2 \,dx,
\end{equation*}
and so we conclude that $\Ateps{u(t)}{v(t)} \le \Ateps{u}{v}$.
\par
For the case $t=0$, by lower semicontinuity we get immediately the result.
\end{proof}

In order to obtain the proof of Theorem \ref{atevol}, we need the following
proposition.

\begin{proposition}
\label{atestbelow}
For all $0 \le s \le t \le 1$, we have that
\begin{eqnarray*}
\Ateps{u(t)}{v(t)}-\Ateps{u(s)}{v(s)} &\ge&
2 \int_{\Om} (\eta_\varepsilon +v^2(t)) \nabla u(t) (\nabla g(t)-\nabla g(s))
\,dx +\\
\nonumber
&& -\sigma(t-s) \int_s^t ||\nabla \dot{g}(\tau)||_{L^2(\Om;\R^N)} \,d\tau
\end{eqnarray*}
where $\sigma$ is an increasing positive function with $\sigma(r) \to 0$ as $r \to 0^+$.
\end{proposition}

\begin{proof}
By Proposition \ref{minimality}, we have
$$
\Ateps{u(s)}{v(s)} \le \Ateps{u(t)-g(t)+g(s)}{v(t)}
$$
so that
\begin{eqnarray*}
\nonumber
\Ateps{u(s)}{v(s)} &\le& \Ateps{u(t)}{v(t)}
-2\int_{\Om}(\eta_\varepsilon +v^2(t)) \nabla u(t)(\nabla g(t)-\nabla g(s))\,dx +\\
&&+ \int_{\Om}(\eta_\varepsilon +v^2(t)) |\nabla g(t)-\nabla g(s)|^2 \,dx.
\end{eqnarray*}
Then we conclude that
\begin{eqnarray}
\label{est2}
\nonumber
\Ateps{u(t)}{v(t)}-\Ateps{u(s)}{v(s)} &\ge&
2\int_{\Om}(\eta_\varepsilon +v^2(t)) \nabla u(t)(\nabla g(t)-\nabla g(s))\,dx +\\
\nonumber
&&- \sigma(t-s) \int_s^t ||\nabla \dot{g}(\tau)||_{L^2} \,d\tau
\end{eqnarray}
where
$$
\sigma(r):= (1+\eta_\varepsilon)
\max_{t-s=r} \int_s^t ||\nabla \dot{g}(\tau)||_{L^2(\Om;\R^N)} \,d\tau,
$$
and so the proof is complete.
\end{proof}

We can now prove Theorem \ref{atevol}.

\begin{proof}[Proof of Theorem \ref{atevol}]
Let us consider the sequence $\delta_n \to 0$ given by Lemma \ref{convv}, and
let us indicate the discrete evolutions $u^{\delta_n}$ and $v^{\delta_n}$
defined in \eqref{discrevol} simply by $u_n$ and $v_n$.
Let us denote also by $u(t)$ and $v(t)$ their limits at time $t$ according to
Lemma \ref{convv} and Lemma \ref{convu}. We have that the maps
$t \to u(t)$ and $t \to v(t)$ are strongly measurable from $[0,1]$ to $H^1(\Om)$;
moreover for all $t \in [0,1]$ we have
$0 \le v(t) \le 1$ in $\Om$, $u(t)=g(t)$, $v(t)=1$ on $\partial_D \Om$ and
$t \to v(t)$ is decreasing from $[0,1]$ to $L^2(\Om)$ so that point $(a)$ is proved.
By construction we get point $(b)$ and by Proposition \ref{minimality}
we get point $(c)$.
\par
Let us pass to condition $(d)$. Let us fix $t \in [0,1]$, and let us divide the
interval $[0,t]$ in $k$ subintervals with endpoints $s_j^k:=\frac{jt}{k}$ where
$j=0,1, \cdots, k$. Let us define $\tilde{u}_k(s):=u(s_{j+1}^k)$, and
$\tilde{v}_k(s):=v(s_{j+1}^k)$ for $s_{j}^k <s \le s_{j+1}^k$. Then, applying
Proposition \ref{atestbelow}, we have
\begin{eqnarray}
\Ateps{u(t)}{v(t)} &\ge& \Ateps{u(0)}{v(0)}+
2 \int_0^t \int_{\Om} (\eta_\varepsilon +\tilde{v}_k^2(\tau))
\nabla\tilde{u}_k(\tau) \nabla \dot{g}(\tau) \,dx \,d\tau +\\
\nonumber
&& -\sigma \left( \frac{t}{k} \right) \int_0^t ||\nabla \dot{g}(\tau)||_{L^2} \,d\tau.
\end{eqnarray}
Since $t \to v(t)$ is monotone decreasing from $[0,1]$ to $L^2(\Om)$, we have that
$\tilde{v}_k(s) \to v(s)$ strongly in $L^2(\Om)$ for a.e. $s \in [0,t]$;
consequently, we have that $\tilde{u}_k(s) \to u(s)$ strongly in $H^1(\Om)$ as
noted in Lemma \ref{convu}.
We conclude by the Dominated Convergence Theorem that
\begin{equation*}
\lim_k  \int_0^t \int_{\Om} (\eta_\varepsilon +\tilde{v}_k^2(\tau))
\nabla\tilde{u}_k(\tau) \nabla \dot{g}(\tau) \,dx \,d\tau =
\int_0^t \int_{\Om} (\eta_\varepsilon +v^2(\tau))
\nabla u(\tau) \nabla \dot{g}(\tau) \,dx \,d\tau.
\end{equation*}
By (\ref{est2}) we deduce that
\begin{equation}
\label{estbelow}
\Ateps{u(t)}{v(t)} \ge \Ateps{u(0)}{v(0)}+
2 \int_0^t \int_{\Om} (\eta_\varepsilon +v^2(\tau))
\nabla u(\tau) \nabla \dot{g}(\tau) \,dx \,d\tau.
\end{equation}
On the other hand, from (\ref{mainestimate}), and since
$\Ateps{u_n(0)}{v_n(0)}=\Ateps{u(0)}{v(0)}$ for all $n$, we deduce
\begin{equation}
\label{estabove}
\limsup_n \Ateps{u_n(t)}{v(t)} \le
\Ateps{u(0)}{v(0)}+
2 \int_0^t \int_{\Om} (\eta_\varepsilon +v^2(\tau))
\nabla u(\tau) \nabla \dot{g}(\tau) \,dx \,d\tau.
\end{equation}
Since by semicontinuity we have for all $t \in [0,1]$
$$
\Ateps{u(t)}{v(t)} \le \liminf_n \Ateps{u_n(t)}{v_n(t)},
$$
by (\ref{estbelow}) and (\ref{estabove}), we conclude that
\begin{equation}
\label{atconvergence}
\lim_n \Ateps{u_n(t)}{v_n(t)}= \Ateps{u(t)}{v(t)}.
\end{equation}
In particular
$$
\Ateps{u(t)}{v(t)}= \Ateps{u(0)}{v(0)}+
2 \int_0^t \int_{\Om} (\eta_\varepsilon +v^2(\tau))
\nabla u(\tau) \nabla \dot{g}(\tau) \,dx \,d\tau,
$$
and this proves point $(d)$.
\end{proof}

\begin{remark}
\label{regularity}
{\rm
The map $\{t \to v(t),\, t \in [0,1]\}$ is decreasing from $[0,1]$ to $L^2(\Om)$, so that
$v$ is continuous with respect to the strong topology of $L^2(\Om)$ at all points
except a countable set. Since
$$
\lambda(t):= \mmeps{v(t)}
$$
is monotone increasing (see Proposition \ref{compactonD}), we conclude that $v:[0,1] \to H^1(\Om)$
is continuous with respect to the strong topology at all points except a countable set.
Then we have $v \in L^\infty([0,1], H^1(\Om))$. Moreover, we have that
$u:[0,1] \to H^1(\Om)$ is continuous at the continuity points of $v$ as observed in Lemma
\ref{convu}. We conclude that $u \in L^\infty ([0,1], H^1(\Om))$.
}
\end{remark}

\begin{remark}
\label{minequil}
{\rm
The minimality property of point $(c)$ of Theorem \ref{atevol} holds indeed in this
stronger form: if $t \in ]0,1]$, for all $(u,v) \in H^1(\Om) \times H^1(\Om)$ with
$0 \le v \le v(s)$ on $\Om$ for all $s<t$, and $u=g(t)$, $v=1$ on $\partial_D \Om$,
we have
$$
\Ateps{u(t)}{v(t)} \le \Ateps{u}{v}.
$$
In fact, if $0 \le v \le v(s)$, by the minimality property of $(u(s),v(s))$ we have
$$
\Ateps{u(s)}{v(s)} \le \Ateps{u+g(s)-g(t)}{v},
$$
so that, letting $s \to t$ and using the continuity of $\Ateps{u(\cdot)}{v(\cdot)}$
we get the result.
\par
This stronger minimality property is the reformulation in the context of the
Ambrosio-Tortorelli functional of the minimality of the cracks required in \cite{FM}
(see the Introduction).
}
\end{remark}

\section{Quasi-static growth of brittle fracture}
\label{evofracture}

In this section, we prove that the evolution for the Ambrosio-Tortorelli
functional converges as $\varepsilon \to 0$ to a quasi-static evolution of brittle
fractures in linearly elastic bodies in the sense of \cite{FL}.
\par
Let $\partial_N \Om:= \partial \Om \setminus \partial_D \Om$ and
let $g \in W^{1,1}([0,1];H^1(\Om))$.
In order to treat in a convenient way the boundary condition as $\varepsilon \to 0$,
let $B$ be an open ball such that $\Omb \subset B$, and let
us pose $\Om':= B \setminus \partial_N \Om$ and $\Om_D:= \Om' \setminus \Omb$.
Let $E$ be an extension operator from $H^1(\Om)$ to $H^1(B)$: we
indicate $Eg(t)$ still by $g(t)$ for all $t \in [0,1]$.
\par
In this enlarged context, the following proposition holds.

\begin{proposition}
\label{atevol2}
There exists a strongly measurable map
$$
\begin{array}{c}
[0,1] \\ t
\end{array}
\begin{array}{c}
\longrightarrow \\ \longmapsto
\end{array}
\begin{array}{c}
H^1(\Om') \times H^1(\Om') \\ (u(t),v(t))
\end{array}
$$
such that for all $t \in [0,1]$ we have
$0 \le v(t) \le 1$ in $\Om'$, $u(t)=g(t)$, $v(t)=1$ on $\Om_D$
and
\begin{itemize}
\item[(a)] for all $0 \le s \le t \le 1: v(t) \le v(s)$;
\item[{}]
\item[{(b)}] for all $(u,v) \in H^1(\Om') \times H^1(\Om')$  with $u=g(0)$,
$v=1$ on $\Om_D$:
\begin{equation}
\label{minimprop0}
\Ateps{u(0)}{v(0)} \le \Ateps{u}{v};
\end{equation}
\item[(c)] for $t \in ]0,1]$ and
for all $(u,v) \in H^1(\Om') \times H^1(\Om')$ with
$0 \le v \le v(t)$ on $\Om'$, and $u=g(t)$, $v=1$ on $\Om_D$:
\begin{equation}
\label{minimprop}
\Ateps{u(t)}{v(t)} \le \Ateps{u}{v};
\end{equation}
\item[{}]
\item[(d)] the function $t \to \Ateps{u(t)}{v(t)}$ is absolutely continuous
and
\begin{equation}
\label{energyeps}
\Ateps{u(t)}{v(t)}=\Ateps{u(0)}{v(0)}+ 2
\int_0^t \int_{\Om'} (\eta_\varepsilon+v^2(\tau))
\nabla u(\tau) \nabla \dot{g}(\tau) \,dx\,d\tau.
\end{equation}
\end{itemize}
\end{proposition}

\begin{proof}
Let us consider the map $t \to (\ueps(t),\veps(t))$ from $[0,1]$
to $H^1(\Om) \times H^1(\Om)$ given by Theorem \ref{atevol}.
Recall that for all $t \in [0,1]$ we have
$\ueps(t)=g(t)$, $v(t)=1$ on
$\partial_D \Om$, and $0 \le \veps(t) \le 1$ in $\Om$.
We extend $\ueps(t)$ and $\veps(t)$ to $\Om'$ posing $\ueps(t)=g(t)$
and $\veps(t)=1$ on $\Om_D$. Then we obtain
a strongly measurable map $t \to (\ueps(t),\veps(t))$ from $[0,1]$
to $H^1(\Om') \times H^1(\Om')$ such that $0 \le \veps(t) \le 1$ in $\Om'$,
$\ueps(t)=g(t)$, $\veps(t)=1$ on $\Om_D$, and such that
$$
\Ateps{\ueps(t)}{\veps(t)} \le \Ateps{u}{v}
$$
for all $(u,v) \in H^1(\Om') \times H^1(\Om')$ with $0 \le v \le \veps(t)$ on
$\Om'$, $u=g(t)$, $v=1$ on $\Om_D$; note in fact that the integrations on $\Om_D$
which appear in both sides are the same. By the same reason, we get
the minimality property at time $t=0$ and deduce that the function
$t \to \Ateps{\ueps(t)}{\veps(t)}$ is absolutely continuous with
$$
\Ateps{\ueps(t)}{\veps(t)}=\Ateps{\ueps(0)}{\veps(0)}+ 2
\int_0^t \int_{\Om'} (\eta_\varepsilon+\veps^2(\tau))
\nabla \ueps(\tau) \nabla \dot{g}(\tau) \,dx\,d\tau.
$$
\end{proof}

From now on, we assume that there exists a constant $C>0$ such that
for all $t \in [0,1]$, $\|g(t)\|_\infty \le C$, and that there exists
$g_h \in W^{1,1}([0,1],H^1(\Om'))$ such that $\|g_h\|_\infty \le C$,
$g_h \in C(\overline{\Om'})$, and
$g_h \to g$ strongly in $W^{1,1}([0,1],H^1(\Om'))$. We indicate by
$(u_{\varepsilon,h},v_{\varepsilon,h})$ the evolution for the Ambrosio-Tortorelli
functional relative to the boundary data $g_h$ given by
Proposition \ref{atevol2}.
The bound on the sup-norm is made in order to apply Ambrosio's compactness theorem
in $SBV$ when $\varepsilon \to 0$.
Notice that we may assume by a truncation argument that
$||\uepsgn(t)||_\infty \le ||g_h(t)||_\infty$, that is
\begin{equation}
\label{linfbound}
||\uepsgn(t)||_\infty \le C.
\end{equation}
We conclude that $\uepsgn(t)$ is
uniformly bounded in $L^\infty(\Om')$ as $\varepsilon$, $h$ and $t$ vary.
Moreover we have that the following holds.

\begin{lemma}
\label{benergy}
There exists a constant $C_1 \ge 0$ depending only on $g$ such that for all
$t \in [0,1]$, $\varepsilon,h$
\begin{equation}
\label{bound}
\Ateps{\uepsgn(t)}{\vepsgn(t)} + \|\uepsgn(t)\|_\infty \le C_1.
\end{equation}
\end{lemma}

\begin{proof}
Notice that $\Ateps{\uepsgn(0)}{\vepsgn(0)} \le
\Ateps{g_h(0)}{1}$ so that
the term $\Ateps{\uepsgn(0)}{\vepsgn(0)}$ is bounded as
$\varepsilon$ and $h$ vary.
We now derive an estimate for the derivative of the total energy.
Since $0 \le \vepsgn(\tau) \le 1$ and $\eta_\varepsilon \to 0$, by
H\"older inequality we get
\begin{multline*}
\left| \int_{\Om'} (\eta_\varepsilon +\vepsgn(\tau)^2)
\nabla \uepsgn(\tau) \nabla \dot{g}_h(\tau) \,dx \right| \le \\
\le
2\left( \int_{\Om'} (\eta_\varepsilon +\vepsgn(\tau)^2) |\nabla \uepsgn(\tau)|^2 \,dx
\right)^\frac{1}{2}
||\nabla \dot{g}_h(\tau)||_{L^2(\Om';\R^N)};
\end{multline*}
since by the minimality property (\ref{minimprop})
$$
\int_{\Om'} (\eta_\varepsilon +\vepsgn(\tau)^2)
|\nabla \uepsgn(\tau)|^2 \,dx \le
\int_{\Om'} (\eta_\varepsilon +\vepsgn(\tau)^2)
|\nabla g_h(\tau)|^2 \,dx,
$$
we get the conclusion by \eqref{energyeps} and \eqref{linfbound}.
\end{proof}

As a consequence of (\ref{bound}), we have
\begin{equation*}
\int_{\Om'} (1-\vepsgn(t))|\nabla \vepsgn(t)| \,dx
\le \mmeps{\vepsgn(t)} \le C_1,
\end{equation*}
so that the functions $w_{\varepsilon,h}(t):= 1-\vepsgn(t)$ have uniformly
bounded variation.
\par
By coarea formula for $BV$-functions (see \cite[Theorem 3.40]{AFP}), we have that
$$
\int_0^1 \hs^{N-1} \left( \partial^* \{\vepsgn(t) >s\} \right) \,ds
=\int_{\Om'} (1-\vepsgn(t))|\nabla \vepsgn(t)| \,dx
$$
($\partial^*$ denotes the essential boundary)
so that by the Mean Value theorem, for all $j \ge 1$ there exists
$b_{\varepsilon,h}^j(t) \in [\frac{1}{2^{j+1}},\frac{1}{2^j}]$ with
\begin{equation}
\label{bper}
\frac{1}{2^{j+1}} \hs^{N-1}
\left( \partial^* \{v_{\varepsilon,n}(t) >b_{\varepsilon,h}^j(t)\}
\right)
\le C_1.
\end{equation}
Let us pose
\begin{equation}
\label{tlevels}
B_{\varepsilon,h}(t):=\left\{b_{\varepsilon,h}^j(t)\,:\,j \ge 1 \right\}.
\end{equation}

We now let $\varepsilon \to 0$. Let $D$ be countable and dense in $[0,1]$
with $0 \in D$.

\begin{lemma}
\label{compact}
There exists a sequence $\varepsilon_n$ such that for all $t \in D$ there
exists $u_h(t) \in SBV(\Om')$, $u_h(t)=g_h(t)$ on $\Om_D$, with
$$
u_{\varepsilon_n,h}(t) \ind{v_{\varepsilon_n,h}(t) >
b_{\varepsilon_n,h}^1(t)} \to u_h(t)
\quad \mbox{ in } SBV(\Om').
$$
In particular for all $t \in D$ we have
\begin{equation}
\label{gradest}
\int_{\Om'} |\nabla u_h(t)|^2 \,dx
+\hs^{N-1}(S_{u_h(t)})+ \|u_h(t)\|_\infty \le C_1.
\end{equation}
\end{lemma}

\begin{proof}
For all $t \in [0,1]$ we may apply Ambrosio's compactness
Theorem \ref{SBVcompact} to the function
$z_n(t):=u_{\varepsilon_n,h}(t) \ind{v_{\varepsilon_n,h}(t) >
b_{\varepsilon_n,h}^1(t)}$: in fact
$z_n(t)$ is bounded in $L^\infty(\Om')$ and
$\nabla z_n(t)$ is bounded in $L^2(\Om')$
by (\ref{bound}), and $S_{z_n(t)} \subseteq
\partial_*\{v_{\varepsilon_n,h}(t) >b_{\varepsilon_n,h}^1(t)\}$ so that
$\hs^{N-1}(S_{z_n(t)})$ is uniformly bounded in $n$ by (\ref{bper}).
Using a diagonal argument, there exists a subsequence
such that for all $t \in D$, $z_n(t) \to u_h(t)$ in
$SBV(\Om')$; in particular, we have that $u_h(t)=g_h(t)$ on $\Om_D$, and by
(\ref{bound}) and the $\Gamma${-}liminf inequality
for the Ambrosio-Tortorelli functional (\ref{atliminf}), we get (\ref{gradest}).
\end{proof}

The following lemma deals with the possibility of truncating at other
levels given by the elements of $B_{\varepsilon_n,h}(t)$.

\begin{lemma}
\label{truncation1}
Let $t \in D$ and $j \ge 1$. For every
$b_{\varepsilon_n,h}^j(t) \in B_{\varepsilon_n,h}(t)$
we have that
$$
u_{\varepsilon_n,h}(t) \ind{v_{\varepsilon_n,h}(t) >
b_{\varepsilon_n,h}^j(t)} \to u_h(t)
\quad {\rm in}\; SBV(\Om').
$$
\end{lemma}

\begin{proof}
Note that, up to a subsequence, $u_{\varepsilon_n,h}(t) \ind{v_{\varepsilon_n,h}(t) >
b_{\varepsilon_n,h}^j(t)} \to z$ in
$SBV(\Om')$ because of Ambrosio's Theorem \ref{SBVcompact}.
By (\ref{linfbound}), we have that
\begin{multline*}
\|u_{\varepsilon_n,h}(t) \ind{v_{\varepsilon_n,h}(t) > b_{\varepsilon_n,h}^j(t)}
-u_{\varepsilon_n,h}(t) \ind{v_{\varepsilon_n,h}(t) > b_{\varepsilon_n,h}^1(t)}\|_{L^2(\Om')} \le \\
\le C \left| \big\{b_{\varepsilon_n,h}^j(t) \le v_{\varepsilon_n,h}(t) \le
b_{\varepsilon_n,h}^1(t)\big\}\right|.
\end{multline*}
Since $v_{\varepsilon_n,h}(t) \to 1$ strongly in $L^2(\Om')$, we conclude that
$$
\left| \big\{b_{\varepsilon_n,h}^j(t) \le v_{\varepsilon_n,h}(t) \le
b_{\varepsilon_n,h}^1(t)\big\} \right| \to 0,
$$
so that
\begin{multline*}
||z-u_h(t)||_{L^2(\Om')}= \\
=\lim_n
||u_{\varepsilon_n,h}(t) \ind{v_{\varepsilon_n,h}(t) > b_{\varepsilon_n,h}^j(t)}
-u_{\varepsilon_n,h}(t) \ind{v_{\varepsilon_n,h}(t) > b_{\varepsilon_n,h}^1(t)}||_{L^2(\Om')}
=0,
\end{multline*}
that is $z=u_h(t)$ and the proof is complete.
\end{proof}

The following lemma deals with the possibility of truncating at time $s$
using the function $v(t)$ for $t \ge s$.

\begin{lemma}
\label{truncation2}
Let $s,t \in D$ with $s \le t$, and $j \ge 1$. Then for every
$b_{\varepsilon_n,h}^j(t) \in B_{\varepsilon_n,h}(t)$ we have that
$$
u_{\varepsilon_n,h}(s) \ind{v_{\varepsilon_n,h}(t) >
b_{\varepsilon_n,h}^j(t)} \to u_h(s)
\quad {\rm in}\; SBV(\Om').
$$
\end{lemma}

\begin{proof}
Up to a subsequence, by Ambrosio's Theorem, we have that
$$
u_{\varepsilon_n,h}(s) \ind{v_{\varepsilon_n,h}(t) > b_{\varepsilon_n,h}^j(t)} \to z
\quad {\rm in}\; SBV(\Om').
$$
Since $v_{\varepsilon_n,h}(t) \le v_{\varepsilon_n,h}(s)$, we have that
$\{v_{\varepsilon_n,h}(t)>b_{\varepsilon_n,h}^j(t)\} \subseteq \{v_{\varepsilon_n,h}(s)>b_{\varepsilon_n,h}^{j+1}(s)\}$.
Then we have
$$
\|u_{\varepsilon_n,h}(s) \ind{v_{\varepsilon_n,h}(t) > b_{\varepsilon_n,h}^j(t)}
-u_{\varepsilon_n,h}(s) \ind{v_{\varepsilon_n,h}(s) > b_{\varepsilon_n,h}^{j+1}(s)}
\|_{L^2(\Om')} \le
C \left|\big\{v_{\varepsilon_n,h}(t) \le b_{\varepsilon_n,h}^j(t)\big\}\right|.
$$
Since $v_{\varepsilon_n,h}(t) \to 1$ strongly in $L^2(\Om')$, we conclude that
$\left|\big\{v_{\varepsilon_n,h}(t) \le b_{\varepsilon_n,h}^j(t)\big\}\right| \to 0$.
By Lemma \ref{truncation1} we have
$$
u_{\varepsilon_n,h}(s) \ind{v_{\varepsilon_n,h}(s) >
b_{\varepsilon_n,h}^{j+1}(s)} \to u_h(s)
\quad {\rm in}\; SBV(\Om'),
$$
so that $z=u_h(s)$ and the proof is complete.
\end{proof}

We now pass to the analysis of $u_h(t)$ with $t \in D$.
The following minimality property for the functions $u_h(t)$
with $t \in D$ is crucial for the subsequent results.

\begin{theorem}
\label{minthm}
Let $t \in D$. Then for every $z \in SBV(\Om')$ with $z=g_h(t)$ on
$\Om_D$, we have that
\begin{equation*}
\int_{\Om'} |\nabla u_h(t)|^2 \,dx \le
\int_{\Om'} |\nabla z|^2 \,dx +
\hs^{N-1} \left( S_z \setminus \bigcup_{s \le t, s \in D} S_{u_h(s)} \right).
\end{equation*}
\end{theorem}

The proof is quite technical, and it is postponed to Section \ref{secminthm}.
We now let $h \to \infty$.

\begin{proposition}
\label{ntoinfty}
There exists $h_n \to \infty$ such that for all $t \in D$ there exists
$u(t) \in SBV(\Om')$ with $u(t)=g(t)$ on $\Om_D$ such that
$u_{h_n}(t) \to u(t)$ in $SBV(\Om')$. Moreover, $\nabla u_{h_n}(t) \to
\nabla u(t)$ strongly in $L^2(\Om';\R^N)$ and for all
$z \in SBV(\Om')$ with $z=g(t)$ on $\Om_D$ we have
\begin{equation*}
\int_{\Om'} |\nabla u(t)|^2 \,dx \le
\int_{\Om'} |\nabla z|^2 \,dx +
\hs^{N-1} \left( S_z \setminus \bigcup_{s \le t, s \in D} S_{u(s)} \right).
\end{equation*}
\end{proposition}

\begin{proof}
The compactness is given by Ambrosio's Theorem in view of \eqref{gradest}.
The strong convergence of the gradients and the minimality
property is a consequence of the minimality property of Theorem
\ref{minthm} and of \cite[Theorem 2.1]{FL}.
\end{proof}

We can now deal with $\varepsilon$ and $h$ at the same time.

\begin{proposition}
\label{compactonD}
There exists $\varepsilon_n \to 0$ and $h_n \to +\infty$ such that for all $t \in D$
there exists $u(t) \in SBV(\Om')$ with $u(t)=g(t)$ on $\Om_D$ such that for all
$j \ge 1$
$$
u_{\varepsilon_n,h_n}(t) 1_{\{v_{\varepsilon_n,h_n}(t) > b^j_{\varepsilon_n,h_n}(t)\}}
\to u(t) \quad \mbox{ in }SBV(\Om').
$$
Furthermore for all $z \in SBV(\Om')$ with $z=g(t)$ on $\Om_D$ we have
\begin{equation*}
\int_{\Om'} |\nabla u(t)|^2 \,dx \le
\int_{\Om'} |\nabla z|^2 \,dx +
\hs^{N-1} \left( S_z \setminus \bigcup_{s \le t, s \in D} S_{u(s)} \right),
\end{equation*}
and we may suppose that the functions $\lambda_{\varepsilon_n,h_n}$
converge pointwise on $[0,1]$ to an increasing function $\lambda$ such that
for all $t \in D$
\begin{equation}
\label{size3}
\lambda(t) \ge \hs^{N-1}
\left( \bigcup_{s \le t, s \in D} S_{u(s)} \right).
\end{equation}
Finally, we have that for all $t \in D$
\begin{equation}
\label{gradestbis}
\int_{\Om'} |\nabla u(t)|^2 \,dx
+\hs^{N-1}(S_{u(t)})+ \|u(t)\|_\infty \le C_1.
\end{equation}
\end{proposition}

\begin{proof}
We find $\varepsilon_n$ and $h_n$ combining Lemma \ref{compact} and
Proposition \ref{ntoinfty}, and using a diagonal argument.
Passing to the second part of the proposition, notice that the
functions $\lambda_{\varepsilon_n,h_n}$ are monotone increasing. In fact
if $s \le t$, since $v_{\varepsilon_n,h_n}(t) \le v_{\varepsilon_n,h_n}(s)$, and
$v_{\varepsilon_n,h_n}(t)=1$ on $\Om_D$, by the minimality property
(\ref{minimprop}), we have that
$$
\Atepsn{u_{\varepsilon_n,h_n}(s)}{v_{\varepsilon_n,h_n}(s)} \le
\Atepsn{u_{\varepsilon_n,h_n}(s)}{v_{\varepsilon_n,h_n}(t)},
$$
so that
\begin{multline*}
\lambda_{\varepsilon_n,h_n}(t)-\lambda_{\varepsilon_n,h_n}(s) \ge \\
\ge \ellepsnp{u_{\varepsilon_n,h_n}(s)}{v_{\varepsilon_n,h_n}(s)}
-\ellepsnp{u_{\varepsilon_n,h_n}(s)}{v_{\varepsilon_n,h_n}(t)} \ge 0.
\end{multline*}
Moreover by (\ref{bound}) we have $0 \le \lambda_{\varepsilon_n,h_n} \le C_1$.
Applying Helly's theorem, we get that there exists an increasing function
$\lambda$ up to a subsequence $\lambda_{\varepsilon_n,h_n} \to \lambda$
pointwise in $[0,1]$. In order to prove \eqref{size3}, let us fix
$s_1, \dots, s_m \in D \cap [0,t]$; we want to prove that
\begin{equation}
\label{sizek}
\lambda(t)=\lim_n \lambda_{\varepsilon_n,h_n}(t) \ge \hs^{N-1}
\left( \bigcup_{i=1}^m S_{u(s_i)} \right).
\end{equation}
Then taking the sup over all possible $s_1, \dots, s_m$, we can deduce
(\ref{size3}). Consider $z_n \in SBV(\Om',\R^m)$ defined as
$$
z_n(x):=(u_{\varepsilon_n,h_n}(s_1), \dots, u_{\varepsilon_n,h_n}(s_m)).
$$
Notice that by (\ref{bound}), and the fact that
$t \to v_{\varepsilon_n,h_n}(t)$ is decreasing
in $L^2(\Om')$, we obtain that there exists $C'>0$ such that for all n
$$
\ellepsnp{z_n(t)}{v_{\varepsilon_n,h_n}(t)} +\mmnp{v_{\varepsilon_n,h_n}(t)} \le C'.
$$
Then we may apply \cite[Lemma 3.2]{Fo} obtaining (\ref{sizek}).
Finally \eqref{gradestbis} is a consequence of \eqref{bound} and the lower
semicontinuity \eqref{atliminf}. The proof is now concluded.
\end{proof}

Let us extend the evolution $\{t \to u(t), t \in D\}$ of Proposition
\ref{compactonD} to the entire interval $[0,1]$. Let us pose for every
$t \in [0,1]$
\begin{equation}
\label{fracture}
\Gamma(t):= \bigcup_{s \in D, s \le t} S_{u(s)}.
\end{equation}

\begin{proposition}
\label{extension}
For every $t \in [0,1]$ there exists $u(t) \in SBV(\Om')$ with $u(t)=g(t)$ on
$\Om_D$ such that $\nabla u \in L^\infty([0,1], L^2(\Om';\R^N))$,
$\nabla u$ is left continuous in $[0,1] \setminus D$ with respect to the
strong topology, and such that, if $\Gamma$ is as in \eqref{fracture},
the following hold:
\begin{itemize}
\item[(a)] for all $t \in [0,1]$
\begin{equation}
\label{jumpt}
S(u(t)) \subseteq
\Gamma(t) \;\mbox{ up to a set of }\hs^{N-1}-\mbox{measure }0,
\end{equation}
and if $\lambda$ is as in Proposition \ref{compactonD}
\begin{equation}
\label{size2}
\lambda(t) \ge \hs^{N-1}(\Gamma(t));
\end{equation}
\item[{}]
\item[(b)] for all $z \in SBV(\Om')$ with $z=g(0)$ on $\Om_D$
\begin{equation}
\label{mint0}
\int_{\Om'} |\nabla u(0)|^2 \,dx + \hs^{N-1} \left( S_u \right) \le
\int_{\Om'} |\nabla z|^2 \,dx +
\hs^{N-1} \left( S_z \right).
\end{equation}
\item[{}]
\item[(c)] for all $t \in ]0,1]$ and for all $z \in SBV(\Om')$ with $z=g(t)$ on
$\Om_D$
\begin{equation}
\label{mint}
\int_{\Om'} |\nabla u(t)|^2 \,dx \le
\int_{\Om'} |\nabla z|^2 \,dx +
\hs^{N-1} \left( S_z \setminus \Gamma(t) \right).
\end{equation}
\end{itemize}
Finally,
\begin{equation}
\label{energybelow}
\Es(t) \ge \Es(0)+
2\int_0^t \int_{\Om'} \nabla u(\tau) \nabla \dot{g}(\tau) \,dx\,d\tau
\end{equation}
where
\begin{equation}
\label{defenergy}
\Es(t):= \int_{\Om'} |\nabla u(t)|^2 \,dx +\hs^{N-1}(\Gamma(t)).
\end{equation}
\end{proposition}

\begin{proof}
Let $t \in [0,1] \setminus D$ and
let $t_n \in D$ with $t_n \nearrow t$; by (\ref{gradestbis}) we can apply
Ambrosio's Theorem obtaining $u \in SBV(\Om')$ with $u=g(t)$
on $\Om_D$ such that $u(t_n) \to u$ in $SBV(\Om')$ up to subsequences.
Let us pose $u(t):=u$. By \cite[Lemma 3.7]{FL}, we have that
\eqref{jumpt} and \eqref{mint} hold, and that
the convergence $\nabla u(t_n) \to \nabla u$ is strong in $L^2(\Om';\R^N)$.
Notice that $\nabla u(t)$ is uniquely determined by (\ref{jumpt}) and (\ref{mint})
since the gradient of the solutions of the minimum problem
$$
\min \left\{ \int_{\Om'} |\nabla u|^2 \,dx \,:\, u=g(t) \;on\;\Om_D,
S_u \subseteq \Gamma(t) \;\mbox{ up to a set of }\hs^{N-1} -\mbox{measure }0\right\}
$$
is unique by the strict convexity of the functional. We conclude that $\nabla u(t)$
is well defined. The argument above proves that $\nabla u$ is left continuous at
all the points of $[0,1] \setminus D$. It turns out that $\nabla u$ is continuous
in $[0,1]$ up to a countable set. In fact let us consider
$t \in [0,1] \setminus (D \cup \ns)$ where $\ns$ is the set of discontinuities
of the function $\hs^{N-1}(\Gamma(\cdot))$. Let $t_n \searrow t$.
By Ambrosio's Theorem, we have that there exists $u \in SBV(\Om')$ with
$u=g(t)$ on $\Om_D$ such that, up to a subsequence, $u(t_n) \to u$ in $SBV(\Om')$.
Since $t$ is a continuity point of $\hs^1(\Gamma(\cdot))$, we deduce that
$S_u \subseteq \Gamma(t)$ up to a set of $\hs^{N-1}$-measure $0$. Moreover
by \cite[Lemma 3.7]{FL} we have that $u$ satisfies the minimality property
(\ref{mint}), and $\nabla u(t_n) \to \nabla u$ strongly in $L^2(\Om';\R^N)$.
We deduce that $\nabla u=\nabla u(t)$, and so $\nabla u(\cdot)$ is continuous in
$[0,1] \setminus (D \cup \ns)$. We conclude that $\nabla u(\cdot)$ is continuous
in $[0,1]$ up to a countable set, so that $\nabla u \in L^\infty([0,1]; L^2(\Om';\R^N))$.
\par
We have that \eqref{size2} is a direct consequence of \eqref{size3}, while
\eqref{mint0} is a consequence of \eqref{minimprop0} and the $\Gamma${-}convergence
result of Ambrosio and Tortorelli \cite{AT1} and \cite{AT2}.
\par
Finally, in order to prove (\ref{energybelow}), we can reason in the following way.
Given $t \in [0,1]$ and $m>0$,
let $s_i^m:= \frac{i}{m}t$ for all $i=0, \ldots, m$. Let us pose
$u^m(s):=u(s_{i+1}^m)$ for $s_i^m< s \le s_{i+1}^m$. By (\ref{mint}) we have
\begin{equation}
\label{energybelowk}
\Es(t) \ge \Es(0)+
2\int_0^t \int_{\Om'} \nabla u^m(\tau) \nabla \dot{g}(\tau) \,d\tau \,dx+o_m
\end{equation}
where $o_m \to 0$ for $m \to +\infty$ because $g$ is absolutely continuous.
Since $\nabla u$ is continuous with respect to the strong topology of
$L^2(\Om';\R^N)$ in $[0,1]$ up to a countable set, passing to the limit for
$m \to +\infty$ we deduce that \eqref{energybelow} holds, and the proof is concluded.
\end{proof}

We are now in position to prove our convergence result.
We need the following lemma.

\begin{lemma}
\label{lemwconvgradt}
Let $\tilde{\ns}$ be the set of discontinuity points of the
function $\lambda$ given by Proposition \ref{compactonD}. Then for every
$t \in [0,1] \setminus \tilde{\ns}$, and $j \ge 1$ we have that
\begin{equation*}
\nabla u_{\varepsilon_n,h_n}(t)
\ind{v_{\varepsilon_n,h_n}(t) > b_{\varepsilon_n,h_n}^j(t)}
\weak \nabla u(t) \quad
{\rm weakly\;in}\; L^2(\Om';\R^N).
\end{equation*}
\end{lemma}

\begin{proof}
Let $t \in [0,1] \setminus \tilde{\ns}$: we may suppose that
$t \notin D$, since otherwise the result has already been established.
Let $s \in D$ with $s <t$. We pose
$$
J:=
\inf \left\{ \int_{\Om'} (\eta_{\varepsilon_n}+v_{\varepsilon_n,h_n}^2(t))
|\nabla z|^2 \,dx
\,:\,z=g_{h_n}(s) \;{\rm on}\;\Om_D \right\},
$$
and we indicate by $w_n(s,t)$ the minimum point of this problem. Notice that
$u_{\varepsilon_n,h_n}(t)-w_n(s,t)$ is the minimum for
$$
K:=
\inf \left\{ \int_{\Om'} (\eta_{\varepsilon_n}+v_{\varepsilon_n,h_n}^2(t))
|\nabla z|^2 \,dx
\,:\,z=g_{h_n}(t)-g_{h_n}(s) \;{\rm on}\;\Om_D \right\}.
$$
Comparing $u_{\varepsilon_n,h_n}(t)-w_n(s,t)$ with $g_{h_n}(t)-g_{h_n}(s)$,
we have
\begin{multline}
\label{compuwt}
\int_{\Om'} (\eta_{\varepsilon_n}+v_{\varepsilon_n,h_n}^2(t))
|\nabla u_{\varepsilon_n,h_n}(t)-\nabla w_n(s,t)|^2 \,dx \le \\
\le \int_{\Om'} (\eta_{\varepsilon_n}+v_{\varepsilon_n,h_n}^2(t))
|\nabla g_{h_n}(t)-\nabla g_{h_n}(s)|^2 \,dx.
\end{multline}
Since $u_{\varepsilon_n,h_n}(s)-w_n(s,t)$ is a good test for $J$, we have
$$
\int_{\Om'} (\eta_{\varepsilon_n}+v_{\varepsilon_n,h_n}^2(t))
\nabla w_n(s,t) (\nabla u_{\varepsilon_n,h_n}(s)-\nabla w_n(s,t)) \,dx=0,
$$
and so the following equality holds
\begin{multline*}
\int_{\Om'} (\eta_{\varepsilon_n}+v_{\varepsilon_n,h_n}^2(t))
(|\nabla u_{\varepsilon_n,h_n}(s)|^2-|\nabla w_n(s,t)|^2) \,dx= \\
= \int_{\Om'} (\eta_{\varepsilon_n}+v_{\varepsilon_n,h_n}^2(t))
(|\nabla u_{\varepsilon_n,h_n}(s)-\nabla w_n(s,t)|^2) \,dx.
\end{multline*}
Since $v_{\varepsilon_n,h_n}(t) \le v_{\varepsilon_n,h_n}(s)$ and by minimality of
$u_{\varepsilon_n,h_n}(s)$ we have
\begin{multline*}
\int_{\Om'} (\eta_{\varepsilon_n}+v_{\varepsilon_n,h_n}^2(t))
|\nabla u_{\varepsilon_n,h_n}(s)|^2 \,dx+\lambda_{\varepsilon_n,h_n}(s) \le \\
\le \int_{\Om'} (\eta_{\varepsilon_n}+v_{\varepsilon_n,h_n}^2(s))
|\nabla u_{\varepsilon_n,h_n}(s)|^2 \,dx+\lambda_{\varepsilon_n,h_n}(s) \le \\
\le
\int_{\Om'} (\eta_{\varepsilon_n}+v_{\varepsilon_n,h_n}^2(t))
|\nabla w_n(s,t)|^2 \,dx+\lambda_{\varepsilon_n,h_n}(t).
\end{multline*}
so that
\begin{multline}
\label{compuws}
\int_{\Om'}
(\eta_{\varepsilon_n}+v_{\varepsilon_n,h_n}^2(t))
(|\nabla u_{\varepsilon_n,h_n}(s)-\nabla w_n(s,t)|^2) \,dx = \\
=\int_{\Om'} (\eta_{\varepsilon_n}+v_{\varepsilon_n,h_n}^2(t))
(|\nabla u_{\varepsilon_n,h_n}(s)|^2-|\nabla w_n(s,t)|^2) \,dx \le \\
\le \lambda_{\varepsilon_n,h_n}(t)-\lambda_{\varepsilon_n,h_n}(s).
\end{multline}
By (\ref{compuwt}) and (\ref{compuws}), we conclude that there exists $C'>0$ with
\begin{multline}
\label{near2}
\int_{\Om'} (\eta_{\varepsilon_n}+v_{\varepsilon_n,h_n}^2(t))
(|\nabla u_{\varepsilon_n,h_n}(t)-\nabla u_{\varepsilon_n,h_n}(s)|^2) \,dx \le \\
\le C'
\|\nabla g_{h_n}(t)-\nabla g_{h_n}(s)\|+
(\lambda_{\varepsilon_n,h_n}(t)-\lambda_{\varepsilon_n,h_n}(s)).
\end{multline}
Then we conclude that for $b_{\varepsilon_n,h_n}^j(t)  \in B_{\varepsilon_n,h_n}(t)$
\begin{multline}
\label{near3}
|| \nabla u_{\varepsilon_n,h_n}(t) \ind{v_{\varepsilon_n,h_n}(t) >
b_{\varepsilon_n,h_n}^j(t)}
-\nabla u_{\varepsilon_n,h_n}(s) \ind{v_{\varepsilon_n,h_n}(t) >
b_{\varepsilon_n,h_n}^j(t)}||_{L^2(\Om';\R^N)} \le \\
\le o(t-s)
\end{multline}
since $\lambda_{\varepsilon_n,h_n} \to \lambda$ pointwise, and
$t$ is a continuity point for $\lambda$.
Recall that by Lemma \ref{truncation2}
$$
\nabla u_{\varepsilon_n,h_n}(s) \ind{v_{\varepsilon_n,h_n}(t) >
b_{\varepsilon_n,h_n}^j(t)} \weak \nabla u(s)
\quad {\rm weakly\;in}\; L^2(\Om';\R^N).
$$
Since
\begin{multline*}
\nabla u_{\varepsilon_n,h_n}(t) \ind{v_{\varepsilon_n,h_n}(t) >
b_{\varepsilon_n,h_n}^j(t)} -\nabla u(t) = \\
=(\nabla u_{\varepsilon_n,h_n}(t) \ind{v_{\varepsilon_n,h_n}(t) >
b_{\varepsilon_n,h_n}^j(t)}
-\nabla u_{\varepsilon_n,h_n}(s) \ind{v_{\varepsilon_n,h_n}(t) >
b_{\varepsilon_n,h_n}^j(t)})+ \\
+(\nabla u_{\varepsilon_n,h_n}(s) \ind{v_{\varepsilon_n,h_n}(t) >
b_{\varepsilon_n,h_n}^j(t)}-\nabla u(s))+(\nabla u(s)-\nabla u(t)),
\end{multline*}
by (\ref{near3}) and the left continuity of $\{\tau \to \nabla u(\tau)\}$ at the points of
$[0,1] \setminus D$, we have that
$$
\nabla u_{\varepsilon_n,h_n}(t) \ind{v_{\varepsilon_n,h_n}(t) >
b_{\varepsilon_n,h_n}^j(t)} \weak \nabla u(t)
\quad {\rm weakly\;in}\; L^2(\Om';\R^N),
$$
so that the lemma is proved.
\end{proof}

We can now pass to the proof of the main theorem of the paper.

\begin{proof}[Proof of Theorem \ref{mainthm}]
By Proposition \ref{atevol2}, we may extend $(\uepsgn(t),\vepsgn(t))$ to $\Om'$
posing $\uepsgn(t)=g_h(t)$ and $\vepsgn(t)=1$ on $\Om_D$, obtaining a quasi-static
evolution in $\Om'$. In this context, the points of $\partial_D \Om$ where the boundary
condition is violated in the limit simply become discontinuity points of the extended
function. Thus we prove the result in this equivalent setting involving $\Om'$.
\par
Let $\varepsilon_n \to 0$ and $h_n \to +\infty$ be the sequences determined
by Proposition \ref{compactonD}. Let us indicate
$u_{\varepsilon_n,h_n}(t), v_{\varepsilon_n,h_n}(t)$ and $F_{\varepsilon_n}$ by
$u_n(t),v_n(t)$ and $F_n$.
Moreover, let us write $B_{n}(t)$ and $b_{n}^j(t)$ for
$B_{\varepsilon_n,h_n}(t)$ and
$b_{\varepsilon_n,h_n}^j(t)$.
Let $\{t \to u(t) \in SBV(\Om')\,,t \in [0,1]\}$ be the evolution
relative to the boundary data $g$ given by Proposition \ref{extension}; up to a
subsequence, we have that $u_n(t) \ind{v_n(t)>b_n^j(t)} \to u(t)$ in $SBV(\Om')$
for all $j \ge 1$ and for all $t$ in a countable and dense subset
$D \subseteq [0,1]$ with $0 \in D$. Moreover we have that
\begin{equation}
\label{energybelow2}
\Es(t) \ge \Es(0)+
2\int_0^t \int_{\Om'} \nabla u(\tau) \nabla \dot{g}(\tau) \,dx\,d\tau,
\end{equation}
where $\Es(t):=\int_{\Om'} |\nabla u(t)|^2 \,dx +\hs^{N-1} \left( \Gamma(t) \right)$
and $\Gamma(t)$ is as in \eqref{fracture}.
\par
By point $(b)$ of Proposition \ref{atevol2} and the Ambrosio-Tortorelli Theorem
we have
\begin{equation}
\label{convat0}
\lim_n F_n(u_n(0),v_n(0)) = \Es(0).
\end{equation}
For $m \ge 1$, notice that
\begin{multline*}
\int_{\Om'} (\eta_{\varepsilon_n}+v_n^2(\tau)) \nabla u_n(\tau) \nabla \dot{g}_{h_n}(\tau) \,dx =
\int_{\Om'} (\eta_{\varepsilon_n}+v_n^2(\tau)) \nabla u_n(\tau)
\ind{v_n(\tau) > b_n^m(\tau)} \nabla \dot{g}_{h_n}(\tau) \,dx +\\
+ \int_{\Om'} (\eta_{\varepsilon_n}+v_n^2(\tau)) \nabla u_n(\tau)
\ind{v_n(\tau) \le b_n^m(\tau)} \nabla \dot{g}_{h_n}(\tau) \,dx.
\end{multline*}
If $\tau \in [0,1]$, we have the estimate
\begin{multline*}
\left| \int_{\Om'} (\eta_{\varepsilon_n}+v_n^2(\tau)) \nabla u_n(\tau)
\ind{v_n(\tau) \le b_k^m(\tau)} \nabla \dot{g}_{h_n}(\tau) \,dx
\right| \le \\
\le \sqrt{\eta_{\varepsilon_n}+ \frac{1}{2^{2m}}}
\left( \int_{\Om'} (\eta_{\varepsilon_n}+v_n^2(\tau)) |\nabla u_n(\tau)|^2 \,dx \right)^{\frac{1}{2}}
||\nabla \dot{g}_{h_n}(\tau)||_{L^2(\Om';\R^N)} \le \\
\le \sqrt{\eta_{\varepsilon_n}+ \frac{1}{2^{2m}}}C \to \frac{C}{2^m}.
\end{multline*}
Moreover, by Lemma \ref{lemwconvgradt} we have that for a.e. $\tau \in [0,1]$
$$
\lim_n \int_{\Om'} (\eta_{\varepsilon_n}+v_n^2(\tau)) \nabla u_n(\tau)
\ind{v_n(\tau) > b_n^m(\tau)} \nabla \dot{g}_{h_n}(\tau) \,dx =
\int_{\Om'} \nabla u(\tau) \nabla \dot{g}(\tau) \,dx,
$$
and we deduce that for such $\tau$
$$
\limsup_n \left|
\int_{\Om'} (\eta_{\varepsilon_n}+v_n^2(\tau)) \nabla u_n(\tau) \nabla \dot{g}_{h_n}(\tau) \,dx -
\int_{\Om'} \nabla u(\tau) \nabla \dot{g}(\tau) \,dx \right| \le \frac{C}{2^m}.
$$
Since $m$ is arbitrary, we have that for a.e. $\tau \in [0,1]$
\begin{equation}
\label{convderiv}
\lim_n
\int_{\Om'} (\eta_{\varepsilon_n}+v_n^2(\tau)) \nabla u_n(\tau) \nabla \dot{g}_{h_n}(\tau) \,dx =
\int_{\Om'} \nabla u(\tau) \nabla \dot{g}(\tau) \,dx.
\end{equation}
By (\ref{energyeps}), (\ref{convat0}), (\ref{convderiv}) and the Dominated
Convergence Theorem, we conclude that for all
$t \in [0,1]$
\begin{equation}
\label{energyconv2}
\lim_n F_n(u_n(t),v_n(t)) =
\Es(0)
+2 \int_0^t \int_{\Om'} \nabla u(\tau) \nabla \dot{g}(\tau) \,dx \,d\tau.
\end{equation}
Since $\liminf_n F_n(u_n(t),v_n(t)) \ge \Es(t)$ by (\ref{atliminf}), by
(\ref{energybelow2}) we have for all $t \in [0,1]$
$$
\lim_n F_n(u_n(t),v_n(t)) = \Es(t).
$$
In particular we get
\begin{equation}
\label{energyeq}
\Es(t)=\Es(0)
+2 \int_0^t \int_{\Om'} \nabla u(\tau) \nabla \dot{g}(\tau) \,dx \,d\tau,
\end{equation}
so that, recalling all the properties stated in Proposition \ref{extension}, we deduce that
$\{t \to u(t)\,:\,t \in [0,1]\}$ is a quasi-static evolution relative to the boundary
data $g$. Point $(a)$ is a consequence of (\ref{energyeq}) and (\ref{energyconv2}).
\par
Let us pass to point $(b)$.
By Lemma \ref{lemwconvgradt}, we know that if $\tilde{\ns}$ is the
set of discontinuity points of $\lambda$, for all $t \in [0,1] \setminus
\tilde{\ns}$ and for all $j \ge 1$ we have
 $\nabla u_n(t) \ind{v_n(t) > b_k^j(t)} \weak \nabla u(t)$
weakly in $L^2(\Om', \R^N)$.
Since
$$
v_n(t) \nabla u_n(t)=
v_n(t) \nabla u_n(t) \ind{v_n(t) > b_n^j(t)}+
v_n(t) \nabla u_n(t) \ind{v_n(t) < b_n^j(t)},
$$
we get immediately that $v_n(t) \nabla u_n(t) \weak \nabla u(t)$ weakly in
$L^2(\Om', \R^N)$.
For all such $t$, we have that
$$
\liminf_n \int_{\Om'} (\eta_{\varepsilon_n}+v^2_n(t)) |\nabla u_n(t)|^2 \,dx \ge \int_{\Om'}
|\nabla u(t)|^2 \,dx,
$$
and by (\ref{size2})
$$
\liminf_n \mmnp{v_n(t)} \ge \hs^{N-1}(\Gamma(t)).
$$
By point $(a)$, we have that the two preceding inequalities are equalities.
In particular, $\lambda$ and $\hs^{N-1}(\Gamma(\cdot))$ coincide up to a
countable set in $[0,1]$. We deduce that $\lambda$ and $\hs^{N-1}(\Gamma(\cdot))$
have the same continuity points, that is $\tilde{\ns}=\ns$.
We conclude that for all $t \in [0,1] \setminus \ns$ we have
$v_n(t) \nabla u_n(t) \to \nabla u(t)$ strongly in
$L^2(\Om', \R^N)$,
$$
\lim_n \int_{\Om'} (\eta_{\varepsilon_n}+v^2_n(t)) |\nabla u_n(t)|^2 \,dx = \int_{\Om'}
|\nabla u(t)|^2 \,dx,
$$
and
$$
\lim_n \mmnp{v_n(t)} = \hs^{N-1}(\Gamma(t)),
$$
so that point $(b)$ is proved, and the proof of the theorem is complete.
\end{proof}

\section{Proof of Theorem \ref{minthm}}
\label{secminthm}

In this section we give the proof of Theorem \ref{minthm} which is an essential
step in the analysis of Section \ref{evofracture}.
For simplicity of notation, for all $t \in D$ we write
$u(t)$, $u_n(t)$ and $v_n(t)$ for $u_h(t)$,
$u_{\varepsilon_n,h}(t)$ and $v_{\varepsilon_n,h}(t)$ respectively.
Moreover, let us write $B_{n}(t)$, $b_{n}^j(t)$ for $B_{\varepsilon_n,h}(t)$ and
$b_{\varepsilon_n,h}^j(t)$, where $B_{\varepsilon_n,h}(t)$ is defined as in
\eqref{tlevels}.
\par
Given $z \in SBV(\Om')$ with $z=g_h(t)$ on $\Om_D$, we want to see that
\begin{equation}
\label{mainminim2}
\int_{\Om'} |\nabla u(t)|^2 \,dx \le
\int_{\Om'} |\nabla z|^2 \,dx +
\hs^{N-1} \left( S_z \setminus \Gamma(t) \right),
\end{equation}
where $g_h(t) \in H^1(\Om') \cap C(\overline{\Om'})$ and
$\Gamma(t)= \bigcup_{s \le t, s\in D} S_{u(s)}$.
\par
The plan is to use the minimality property (\ref{minimprop}) of the approximating
evolution, so that the main point is to construct a sequence
$(z_n,v_n) \in H^1(\Om') \times H^1(\Om')$ such that $z_n=g_h(t)$, $v_n=1$ on $\Om_D$,
$0 \le v_n \le v_n(t)$, and such that
$$
\lim_n \int_{\Om'} (\eta_n+v_n^2) |\nabla z_n|^2 \,dx=
\int_{\Om'} |\nabla z|^2 \,dx
$$
and
$$
\limsup_n
\left[
MM_n(v_n)-MM_n(v_n(t))
\right]
\le \hs^{N-1} \left( S_z \setminus \Gamma(t) \right),
$$
where we use the notation
$$
MM_n(w):= \mmnp{w}.
$$
If a sequence with these properties exists, then by property (\ref{minimprop})
we get the result.
The following lemma contains the main ideas in order to prove Theorem \ref{minthm}.

\begin{lemma}
\label{minlem}
Let $t \in D$; given $z \in SBV(\Om')$ with $z=g_h(t)$ on $\Om_D$ we have that
\begin{equation}
\label{mineq}
\int_{\Om'} |\nabla u(t)|^2 \,dx \le
\int_{\Om'} |\nabla z|^2 \,dx +
\hs^{N-1} \left( S_z \setminus S_{u(t)} \right).
\end{equation}
\end{lemma}

In order to prove Lemma \ref{minlem}, we need several preliminary results.
Let $z \in SBV(\Om')$ be such that $z=g_h(t)$ on $\Om_D$.
Given $\sigma>0$, let $U$ be a neighborhood of $S_{u(t)}$ such that
$|U| \le \sigma$, and $||\nabla z||_{L^2(U;\R^N)} \le \sigma$.
Let $C:=\{x \in \partial_D \Om: \partial_D \Om \;
{\rm is \; not\; differentiable\; at}\; x\}$.
We recall that there exists a countable and dense set $A \subseteq \R$
such that up to a set of $\hs^{N-1}$-measure zero
$$
S_{u(t)}= \bigcup_{a,b \in A} \partial^* E_a \cap \partial^*E_b
$$
where $E_a:=\{x \in \Om'\,:\, u(t)(x) \ge a\}$ and $\partial^*$ denotes the
essential boundary.
Consider
$$
J_j:= \left\{x \in S_{u(t)} \setminus C: [u(t)(x)] \ge \frac{1}{j} \right\},
$$
with $j$ chosen in such a way that $\hs^{N-1}(S_{u(t)} \setminus J_j) \le \sigma$.
For $x \in J_j$, let $a_1(x), a_2(x) \in A$ be such that
$u^-(t)(x) < a_1(x) < a_2(x) < u^+(t)(x)$ and $a_2(x)-a_1(x) \ge \frac{1}{2j}$.
Following \cite[Theorem 2.1]{FL}, we consider a finite disjoint collection of closed
cubes $\{Q_i\}_{i=1, \ldots,k}$ with center $x_i \in J_j$, radius $r_i$ and with normal
$\nu(x_i)$ such that $\bigcup_{i=1}^k Q_i \subseteq U$,
$\hs^{N-1}(J_j \setminus \bigcup_{i=1}^k Q_i) \le \sigma$, and for
all $i=1, \ldots,k$, $j=1,2$
\begin{itemize}
\item[1.] $\hs^{N-1} \left( S_{u(t)} \cap \partial Q_i \right)=0$;
\item[{}]
\item[2.] $r_i^{N-1} \le c \hs^{N-1} \left( S_{u(t)} \cap Q_i \right)$ for some constant $c>0$;
\item[{}]
\item[3.] $\hs^{N-1} \left( \big[ S_{u(t)} \setminus \partial^* E_{a_j(x_i)} \big] \cap Q_i \right)
\le \sigma r_i^{N-1}$;
\item[{}]
\item[4.] $\hs^{N-1} \left( \big\{y \in \partial^* E_{a_j(x_i)} \cap Q_i\,:\,
{\rm dist}(y,H_i) \ge \frac{\sigma}{2} r_i \big\} \right) <\sigma r_i^{N-1}$ where
$H_i$ denotes the intersection of $Q_i$ with the hyperplane
through $x_i$ orthogonal to $\nu(x_i)$;
\item[{}]
\item[5.] $\hs^{N-1} \left( \big( S_z \setminus S_{u(t)} \big) \cap Q_i \right) < \sigma r_i^{N-1}$ and
$\hs^{N-1} (S_z  \cap \partial Q_i)=0$.
\end{itemize}
Note that we may suppose that $Q_i \subseteq \Om$ if $x_i \in \Om$.
Moreover we may require that (see \cite[Theorem 2.1]{FL} and references therein) for all $i=1, \ldots,k$
and $j=1,2$
\begin{equation}
\label{density}
\| 1_{E_{a_j(x_i)} \cap Q_i} -1_{Q^+_i}\|_{L^1(\Om')}
\le \sigma^2 r_i^N.
\end{equation}
Let us indicate by $R_i$ the rectangle given by the intersection of $Q_i$ with
the strip centered at $H_i$ with width $2\sigma r_i$, and let
us pose $V_i:=\{y+s \nu(x_i):\, y\in \partial Q_i, s \in \R\} \cap R_i$.
Note that up to changing the strip, we can suppose
$\hs^{N-1}(\partial R_i \cap (S_u \cup S_z))=\emptyset$.
\par
If $x_i \in \partial_D \Om$, since $x_i \not\in C$, we may require that
\begin{equation}
\label{distbdry}
{\rm dist}(\partial \Om \cap Q_i, H_i) < \sigma r_i;
\end{equation}
moreover, if $(Q_i^+ \setminus R_i) \subseteq \Om$, we can assume that
$g_h(t) <a_1(x_i)$ on $\partial \Om \cap Q_i$ because $g_h(t)$ is continuous
and $g_h(t)(x_i)=u^-(x_i) <a_1(x_i)$. Similarly we may require that
$g_h(t) >a_2(x_i)$ on $\partial \Om \cap Q_i$ in the case
$(Q_i^- \setminus R_i) \subseteq \Om$.
\par
Since we can reason up to subsequences of $\varepsilon_n$, we may suppose
that $\sum_n \varepsilon_n \le \frac{1}{8}$.
Since by~\eqref{bound} we have that $||u_n(t)||_\infty < C_1$ and $v_n(t) \to 1$
strongly in $L^2(\Om')$,
by Lemma \ref{truncation1} we deduce that $u_n(t) \to u(t)$ in measure.
By \eqref{density}, we deduce that for $n$ large enough
\begin{equation}
\label{densityn}
|Q_i^+ \setminus E^n_{a_2(x_i)}| \le 2\sigma^2 r_i^N,
\end{equation}
where we use the notation $E^n_a:=\{x \in \Om'\,:\,u_n(t)(x) \ge a\}$.
Let $G_n \subseteq ]\frac{\sigma}{4}r_i, \frac{\sigma}{2}r_i[$ be the set of all
$s$ such that
$$
\int_{H_i(s)} (\eta_n+ v_n^2(t)) |\nabla u_n(t)|^2 \,d\hs^{N-1} \ge
\frac{C_1}{\sigma r_i \varepsilon_n};
$$
we get immediately by \eqref{bound} that
$$
|G_n| \le \sigma r_i \varepsilon_n,
$$
so that, posing $G:= \bigcup_n G_n$, we have $|G| \le \frac{\sigma}{8}r_i$ and
$|\,]\frac{\sigma}{4}r_i, \frac{\sigma}{2} r_i[ \,\setminus G\,| \ge \frac{\sigma}{8}r_i$.
From (\ref{densityn}), applying Fubini's Theorem we obtain
$$
\int_{]\frac{\sigma}{4}r_i, \frac{\sigma}{2} r_i[ \setminus G}
\hs^{N-1} \left( H_i(s) \setminus E^n_{a_2(x_i)} \right) ds \le 2\sigma^2 r_i^N,
$$
so that there exists $\overline{s} \in ]\frac{\sigma}{4}r_i,\frac{\sigma}{2} r_i[\,
\setminus G$ such that, posing $H_i^+:=H_i(\overline{s})$, we have
\begin{equation}
\label{h+est}
\hs^{N-1} \left( H_i^+ \setminus E^n_{a_2(x_i) -\frac{\delta}{2}} \right) \le
16\sigma r_i^{N-1}.
\end{equation}
Moreover we have by construction
\begin{equation}
\label{gradh+}
\int_{H_i^+} (\eta_n+v_n^2(t))|\nabla u_n|^2 \,d\hs^{N-1} \le K_n,
\end{equation}
where $K_n$ is of the order of $\frac{1}{\varepsilon_n}$.
In a similar way, there exists $H^-_i:=H_i(\tilde{s})$ with
$\tilde{s} \in ]-\frac{\sigma}{2} r_i, -\frac{\sigma}{4}r_i[$ and
\begin{equation}
\label{h-est}
\hs^{N-1} \left( H_i^- \cap E^n_{a_1(x_i) +\frac{\delta}{2}} \right) \le 16\sigma r_i^{N-1},
\end{equation}
and
\begin{equation}
\label{gradh-}
\int_{H_i^-} (\eta_n+v_n^2(t))|\nabla u_n|^2 \,d\hs^{N-1} \le
K_n
\end{equation}
where $K_n$ is of the order of $\frac{1}{\varepsilon_n}$.
We indicate by $\tilde{R}_i$ the intersection of $Q_i$ with the strip determined
by $H_i^+$ and $H_i^-$.
\par
A similar argument prove that, up to reducing $Q_i$ (preserving the estimates previously
stated), we may suppose that
\begin{equation}
\label{boundvu}
\int_{V_i} (\eta_n+ v_n^2(t)) |\nabla u_n(t)|^2 \,d\hs^{N-1} \le K_n,
\end{equation}
where $K_n$ is of the order of $\frac{1}{\varepsilon_n}$.
\par
In order to prove Lemma \ref{minlem}, we claim that we can suppose $z=g_h(t)$ on
$\Om_D$ and in a neighborhood $\mathcal V$ of
$\partial_D \Om \setminus \bigcup_{i=1}^k Q_i$, $S_z \setminus \bigcup_{i=1}^k R_i$
polyhedral with closure contained in $\Om$, and
$\hs^{N-1}((S_z \setminus S_{u(t)}) \cap Q_i) \le \sigma r_i^{N-1}$ for all
$i=1, \ldots,k$. In fact, by Proposition \ref{regularization}, there exists
$w_m \in SBV(\Om')$ with $w_m=g_h(t)$ in $\Om' \setminus \Omb$ and in a
neighborhood $\mathcal{V}_m$ of $\partial_D \Om$ such that $w_m \to z$
strongly in $L^2(\Om')$, $\nabla w_m \to \nabla z$ strongly in
$L^2(\Om';\R^N)$, $\overline{S_{w_m}} \subseteq \Om$ polyhedral, and such that
for all $A$ open subset of $\Om'$ with $\hs^{N-1}(\partial A \cap S_{z})=0$,
we have
$$
\lim_m \hs^{N-1}(A \cap S_{w_m})=\hs^{N-1}(A \cap S_{z}).
$$
Let us fix $\sigma'>0$ and let us consider for all $i=1, \ldots, k$ a rectangle
$R'_i$ centered in $x_i$, oriented as $R_i$ and such that
$R'_i \subseteq \tint{R_i}$, $\hs^{N-1}(\partial R'_i \cap S_{z})=0$,
$\hs^{N-1}(S_z \cap (\tint{R_i} \setminus R'_i)) \le \sigma' r_i^{N-1}$,
where $\tint{R_i}$ denotes the interior part of $R_i$.
Let $\psi_i$ be a smooth function such that $0 \le \psi_i \le 1$, $\psi_i=1$ on
$R'_i$ and $\psi_i=0$ outside $R_i$. Posing $\psi:=\sum_{i=1}^k \psi_i$, let
us consider $z_m:= \psi z+ (1-\psi) w_m$. Note that
$z_m \to z$ strongly in $L^2(\Om')$, $\nabla z_m \to \nabla z$ strongly in
$L^2(\Om;\R^N)$, $z_m=g_h(t)$ in $\Om_D$ and in a
neighborhood $\mathcal{V}'_m$ of $\partial_D \Om \setminus \bigcup_{i=1}^k R_i$,
$S_{w_m} \setminus \bigcup_{i=1}^k R_i$ is polyhedral with closure contained in $\Om$.
Finally, for $m \to +\infty$, we have
$\hs^{N-1}(S_{z_m} \setminus \bigcup_{i=1}^k Q_i) \to
\hs^{N-1}(S_{z} \setminus \bigcup_{i=1}^k Q_i))$ and
$\limsup_m \hs^{N-1}(S_{z_m} \cap (\tint{R_i} \setminus R_i'))
\le 2\hs^{N-1}(S_z \cap (\tint{R_i} \setminus R_i')) \le 2\sigma' r_i^{N-1}$. So, if
\eqref{mineq} holds for $z_m$, we obtain for $m \to +\infty$
that~\eqref{mineq} holds also for $z$ since $\sigma'$ is arbitrary, and so
the claim is proved.
\par
We begin with the following lemma.

\begin{lemma}
\label{wn}
Let $B_n(t)$ be as in \eqref{tlevels}, and
let us consider $b_n^2:=b_n^{j_2}(t),
b_n^3:=b_n^{j_3}(t) \in B_n(t)$ with $j_2>j_3>1$.
Suppose that $k_n:= \frac{b_n^3}{b_n^2}>1$ and let $k,b$ be such that
$1< k \le k_n$, $b_n^3 \le b$ for all $n$.
Then posing
$$
w_n:=
\left\{
\begin{array}{l}
\frac{k_n}{k_n-1} (v_n(t)-b_n^{3})+b_n^{3} \\ \\
0 \\ \\
v_n(t)
\end{array}
\begin{array}{l}
{\rm in}\; \{b_n^{2} \le v_n(t) \le b_n^{3}\} \\ \\
{\rm in}\; \{v_n(t) \le b_n^{2}\} \\ \\
{\rm in}\; \{v_n(t) \ge b_n^{3}\}
\end{array}
\right.
$$
we have that $w_n \in H^1(\Om')$ with $w_n=1$ on $\Om_D$,
$0 \le w_n \le v_n(t)$ in $\Om'$ and
\begin{equation}
\label{mm1}
\limsup_n \left(
MM_n(w_n)-MM_n(v_n(t)) \right) \le
\frac{2C_1k}{(k-1)^2}+\frac{C_1}{(k-1)(1-b)^2}+\frac{C_1b}{(1-b)^2},
\end{equation}
where $C_1$ is given by (\ref{bound}). Moreover there exist
$b_n^1:=b_n^{j_1}(t) \in B_n(t)$ with $j_1>j_2+1$ and a cut-off function
$\varphi_n \in H^1(\Om')$ with $\varphi_n=0$ in $\{v_n(t) \le b_n^1\}$,
$\varphi_n=1$ on $\{v_n(t) \ge b_n^2\}$ (in particular on $\Om_D$) and such that
\begin{equation}
\label{grad1}
\lim_n \eta_n \int_{\Om'} |\nabla \varphi_n|^2 \,dx = 0
\end{equation}
\end{lemma}

\begin{proof}
$w_n$ is well defined in $H^1(\Om')$, and by construction $w_n=1$ on $\Om_D$ and
$0 \le w_n \le v_n(t)$ in $\Om'$. Let us estimate $MM_n(w_n) -MM_n(v_n)$. Since
$$
\frac{\varepsilon_n}{2} \int_{\Om'} |\nabla w_n|^2 \,dx=
\frac{\varepsilon_n}{2} \int_{\{v_n(t) \ge b_n^{3}\}} |\nabla v_n(t)|^2 \,dx
+
\frac{\varepsilon_n}{2} \int_{\{b_n^{2} \le v_n(t) \le b_n^{3}\}}
|\nabla w_n|^2 \,dx,
$$
and $MM_n(v_n(t)) \le C_1$ by \eqref{bound}, we have that
\begin{eqnarray*}
&{}&
\frac{\varepsilon_n}{2}
\int_{\Om'} |\nabla w_n|^2 \,dx -
\frac{\varepsilon_n}{2} \int_{\Om'} |\nabla v_n(t)|^2 \,dx \le \\
&\le&
\frac{\varepsilon_n}{2} \int_{\{b_n^{2} \le v_n(t) \le b_n^{3}\}}
\left( \frac{k_n^2}{(k_n-1)^2}-1 \right) |\nabla v_n(t)|^2 \,dx
-\frac{\varepsilon_n}{2} \int_{\{v_n(t) \le b_n^{2}\}}
|\nabla v_n(t)|^2 \,dx \le \\
&\le&
C_1\left( \frac{k_n^2}{(k_n-1)^2}-1 \right)=
\frac{C_1(2k_n-1)}{(k_n-1)^2} \le
\frac{2C_1k}{(k-1)^2}.
\end{eqnarray*}
Moreover we have that
\begin{eqnarray*}
&\displaystyle \frac{1}{2\varepsilon_n}& \int_{\Om'} (1-w_n)^2 \,dx -
\frac{1}{2\varepsilon_n} \int_{\Om'} (1-v_n(t))^2 \,dx = \\
&=&
\frac{1}{2\varepsilon_n} \int_{\Om'}
\left[ (1-w_n)^2-(1-v_n(t))^2 \right] \,dx =\\
&=&
\frac{1}{2\varepsilon_n} \int_{\Om'} (v_n(t)-w_n)(2-v_n(t)-w_n) \,dx =\\
&=&
\frac{1}{2\varepsilon_n} \int_{\{b_n^{2} \le v_n(t) \le b_n^{3}\}}
\left( v_n(t)- \frac{k_n}{k_n-1}(v_n(t)-b_n^3)-b_n^{3} \right)(2-v_n(t)-w_n) \,dx +\\
&&
+\frac{1}{2\varepsilon_n} \int_{\{v_n(t) \le b_n^{2}\}} v_n(t)(2-v_n(t))\,dx = \\
&=&
\frac{1}{2\varepsilon_n}\int_{\{b_n^{2} \le v_n(t) \le b_n^{3}\}}
\frac{1}{k_n-1} (b_n^{3}-v_n(t))(2-v_n(t)-w_n)\,dx + \\
&&
+\frac{1}{2\varepsilon_n} \int_{\{v_n(t) \le b_n^{2}\}} v_n(t)(2-v_n(t))\,dx \le \\
&\le&
\frac{C_1}{(k_n-1)(1-b_n^3)^2}+\frac{C_1b_n^2}{(1-b_n^2)^2} \le
\frac{C_1}{(k-1)(1-b)^2}+\frac{C_1b}{(1-b)^2}
\end{eqnarray*}
because $\frac{|\{v_n(t) \le s\}|}{\varepsilon_n} \le \frac{C_1}{(1-s)^2}$.
We conclude that
\begin{equation*}
\limsup_n \left(
MM_n(w_n)-MM_n(v_n(t)) \right) \le
\frac{2C_1k}{(k-1)^2}+\frac{C_1}{(k-1)(1-b)^2}+\frac{C_1b}{(1-b)^2}.
\end{equation*}
Let $j_1>j_2+1$: we have that $b_n^1:=b_n^{j_1}$ and $b_n^2$ are not in adjacent
intervals, and so there exists $l>0$ with $0<l \le b_n^2-b_n^1$.
Let us divide the interval $[b_n^{1},b_n^{2}]$ in
$h_n$ intervals of the same size $I_j, j=1, \dots,h_n$, with $h_n$
such that $\frac{\eta_n}{\varepsilon_n} h_n \to 0$. Since
$$
\sum_{j=1}^{h_n}
\frac{\varepsilon_n}{2} \int_{\{v_n(t) \in I_j\}} |\nabla v_n(t)|^2 \,dx \le
\frac{\varepsilon_n}{2} \int_{\Om'} |\nabla v_n(t)|^2 \,dx \le C_1,
$$
we deduce that there exists $I_n$ such that
\begin{equation}
\label{estcutoff1}
\frac{\varepsilon_n}{2} \int_{\{v_n(t)\in I_n\}}
|\nabla v_n(t)|^2 \,dx \le \frac{C_1}{h_n}.
\end{equation}
Let $\alpha_n, \beta_n$ be the extremes of $I_n$. Let us pose
\begin{equation}
\label{cutoff}
\varphi_n:= \frac{1}{\beta_n-\alpha_n} (v_n-\alpha_n)^+ \wedge 1.
\end{equation}
Then $\varphi_n \in H^1(\Om)$, $\varphi_n=0$ in $\{v_n(t) \le b_n^1\}$,
$\varphi_n=1$ on $\{v_n(t) \ge b_n^2\}$ (in particular on $\Om_D$) and
by \eqref{estcutoff1} and the choice of $h_n$ we have that
\begin{multline*}
\label{grad1'}
\eta_n \int_{\Om'} |\nabla \varphi_n|^2 \,dx =
\eta_n \int_{\{\alpha_n \le v_n(t) \le \beta_n\}} \frac{1}{(\beta_n-\alpha_n)^2}
|\nabla v_n(t)|^2 \,dx
\le \frac{\eta_n}{\varepsilon_n} \frac{2C_1}{h_n} \frac{h_n^2}{l^2} \to 0,
\end{multline*}
so that the proof is complete.
\end{proof}

Let $b_n^1$ be as in Lemma \ref{wn} and let $\delta:=\frac{1}{8j}$ so
that for all $i=1,\ldots,k$
\begin{equation*}
a_1(x_i) < a_1(x_i)+\delta < a_2(x_i)-\delta < a_2(x_i).
\end{equation*}

\begin{lemma}
\label{wni2}
For each $i=1, \ldots, k$, there exists $w_n^{2,i} \in H^1(Q_i)$ and
$[\gamma_n^i-\tau_n^i,\gamma_n^i+\tau_n^i]
\subseteq [a_1(x_i)+\delta, a_2(x_i)-\delta]$ such that
$0 \le w_n^{2,i}\le 1$, $w_n^{2,i}=0$ in
$\{\gamma_n^i-\tau_n^i \le  u_n(t) \le \gamma_n^i+\tau_n^i\} \cap Q_i$, $w_n^{2,i}=1$
on $\left[ \big\{ u_n(t) \le a_1(x_i)+\frac{3}{4}\delta \big\} \cup
\big\{ u_n(t) \ge a_2(x_i)- \frac{3}{4}\delta \big\}\right] \cap Q_i$,
and
\begin{equation}
\label{mm2}
\limsup_n \sum_{i=1}^k MM_n(w_n^{2,i})_{|\{v_n(t) > b_n^1\}} \le o(\sigma).
\end{equation}
Moreover there exists $\varphi_n^{2,i} \in H^1(Q_i)$ such that
$0 \le \varphi_n^{2,i} \le 1$, $\varphi_n^{2,i}=0$ on
$\{\gamma_n^i-\frac{\tau_n^i}{2} \le u_n(t) \le \gamma_n^i+\frac{\tau_n^i}{2}\}
\cap Q_i$,
$\varphi_n^{2,i}=1$ on
$\left[\{u_n(t) \le \gamma_n^i-\tau_n^i\} \cup \{u_n(t) \ge \gamma_n^i+\tau_n^i\}
\right] \cap Q_i$, and
\begin{equation}
\label{grad2}
\lim_n \eta_n
\int_{Q_i \cap \{v_n(t) > b_n^1\}} |\nabla \varphi_n^{2,i}|^2 \,dx = 0.
\end{equation}
\end{lemma}

\begin{proof}
For each $i$ let us consider the strip
$$
S_n^i:=E^n_{a_1(x_i)+\delta} \setminus E^n_{a_2(x_i)-\delta}.
$$
Let $h_n \in \N$ and let us divide
$[a_1(x_i)+\delta, a_2(x_i)-\delta]$
in $h_n$ intervals of the same size: there exists a subinterval with extremes $\alpha^i_n$ and $\beta^i_n$ such that,
posing $\tilde{S}^i_n:=\{x \in \Om'\,: \alpha_n^i \le u_n(t) \le \beta_n^i\}$,
\begin{align}
\label{trick}
&
\int_{\tilde{S}^i_n \cap Q_i}
\left[
\sigma  (\eta_n+v_n^2(t)) |\nabla u_n(t)|^2 +(1-\sigma)
\right]  \,dx \le \\
\nonumber
&\le
\frac{1}{h_n}
\int_{S^i_n \cap Q_i}
\left[
\sigma  (\eta_n+v_n^2(t)) |\nabla u_n(t)|^2 +(1-\sigma)
\right]  \,dx.
\end{align}
Let $\gamma_n^i:= \frac{\alpha_n^i+\beta_n^i}{2}$ and
$\tau^i_n:= \frac{a_2(x_i)-a_1(x_i)-2\delta}{4h_n}$. We pose
\begin{equation*}
w^{2,i}_n:=
\left\{
\begin{array}{l}
\frac{1}{\beta_n^i-\gamma_n^i-\tau_n^i} (u_n(t)-\gamma_n^i-\tau_n^i)^+ \wedge 1 \\ \\
0 \\ \\
\frac{1}{\gamma_n^i-\tau_n^i-\alpha_n^i} (u_n(t)-\gamma_n^i+\tau_n^i)^- \wedge 1
\end{array}
\begin{array}{l}
{\rm in}\; \{u_n(t) \ge \gamma_n^i+\tau_n^i\} \cap Q_i \\ \\
{\rm in}\; \{\gamma_n^i-\tau_n^i \le u_n(t) \le \gamma_n^i+\tau_n^i\} \cap Q_i \\ \\
{\rm in}\; \{u_n(t) \le \gamma_n^i-\tau_n^i\} \cap Q_i.
\end{array}
\right.
\end{equation*}
We have that
\begin{eqnarray*}
\frac{\varepsilon_n}{2}
\int_{Q_i \cap \{v_n(t) > b_n^1\}} |\nabla w_n^{2,i}|^2 \,dx
+
\frac{1}{2\varepsilon_n}
\int_{Q_i \cap \{v_n(t) > b_n^1\}} (1-w_n^{2,i})^2 \,dx\le \\
\le \frac{\varepsilon_n}{2}
\left(
\frac{4h_n^2}{\delta^2}
\int_{\tilde{S}^i_n \cap(Q_i \cap \{v_n(t) > b_n^1\})} |\nabla u_n(t)|^2 \,dx
\right)
+
\frac{1}{2\varepsilon_n}
|\tilde{S}^i_n \cap(Q_i \cap \{v_n(t) > b_n^1\})|.
\end{eqnarray*}
Since by (\ref{trick})
$$
\int_{\tilde{S}^i_n \cap Q_i}
(\eta_n+v_n^2(t)) |\nabla u_n(t)|^2 \,dx \le
\frac{1}{h_n}
\left[
\int_{S^i_n \cap Q_i}
(\eta_n+v_n^2(t)) |\nabla u_n(t)|^2 \,dx
+\frac{1-\sigma}{\sigma} |S^i_n \cap Q_i|
\right]
$$
and
$$
|\tilde{S}^i_n \cap Q_i| \le \frac{1}{h_n}
\left[
\frac{\sigma}{1-\sigma}
\int_{S^i_n \cap Q_i}
(\eta_n+v_n^2(t)) |\nabla u_n(t)|^2 \,dx
+|S^i_n \cap Q_i|
\right]
$$
we have
\begin{multline*}
 MM_n(w_n^{2,i})_{|\{v_n(t) > b_n^1\}} \le \\
\le
\frac{2 h_n \varepsilon_n}{\delta^2 (\eta_n+(b_n^1)^2)}
\left[
\int_{S^i_n \cap Q_i}
(\eta_n+v_n^2(t)) |\nabla u_n(t)|^2 \,dx
+\frac{1-\sigma}{\sigma} |S^i_n \cap Q_i|
\right] +\\
+
\frac{1}{2\varepsilon_n h_n}
\left[
\frac{\sigma}{1-\sigma}
\int_{S^i_n \cap Q_i}
(\eta_n+v_n^2(t)) |\nabla u_n(t)|^2 \,dx
+|S^i_n \cap Q_i|
\right].
\end{multline*}
Summing on $i=1,\ldots,k$, recalling (\ref{bound}) and letting $d \in ]0,1]$ with
$\eta_n +(b_n^1)^2 \ge d^2$ for all $n$, we obtain
\begin{align*}
\sum_{i=1}^k & MM_n(w_n^{2,i})_{|\{v_n(t) > b_n^1\}} \le \\
& \le h_n \varepsilon_n\frac{2}{\delta^2  d^2}
\left[
C_1+\frac{1-\sigma}{\sigma} |\cup Q_i|
\right] +
\frac{1}{\varepsilon_n h_n}
\frac{1}{2}
\left[
\frac{\sigma}{1-\sigma}C_1
+|\cup Q_i|
\right].
\end{align*}
We choose $h_n$ in such a way that the preceding quantity is less than
(recall that $|\cup Q_i| \le |U| <\sigma$)
$$
\sqrt{
\frac{1}{\delta^2 d^2}
\left( C_1+1-\sigma \right) \left( \frac{\sigma}{1-\sigma}C_1+\sigma \right)}.
$$
Then we obtain
\begin{equation*}
\sum_{i=1}^k MM_n(w_n^{2,i})_{|\{v_n(t) > b_n^1\}} \le
\sqrt{
\frac{1}{\delta^2 d^2}
(C_1+1-\sigma))
\left( \frac{\sigma}{1-\sigma}C_1+\sigma \right)}  =o(\sigma).
\end{equation*}
This prove the first part of the lemma.
\par
Let us define $\varphi_n^{2,i}$ as $w_n^{2,i}$ but operating with the levels
$\gamma_n^i-\tau_n^i \le \gamma_n^i-\frac{\tau_n^i}{2}$ and
$\gamma_n^i+\frac{\tau_n^i}{2} \le \gamma_n^i+\tau_n^i$.
Reasoning as above we obtain
\begin{equation*}
\eta_n \sum_{i=1}^k
\int_{Q_i \cap \{v_n(t) > b_n^1\}} |\nabla \varphi_n^{2,i}|^2 \,dx
\le
\frac{16\eta_n h_n}{\delta^2 d^2} (C_1+1-\sigma) \to 0
\end{equation*}
since $h_n$ has been chosen of the order of $\frac{1}{\varepsilon_n}$.
\end{proof}

\begin{lemma}
\label{wn3}
Let $Q_i \subseteq \Om$. Then there exists $w_n^{3,i} \in H^1(Q_i)$ such that
$0 \le w_n^{3,i} \le 1$, $w_n^{3,i}=0$ in a neighborhood of
$H_i^+ \setminus E_{a_2(x_i)-\frac{3}{4}\delta}^n$ and of
$H_i^- \cap E_{a_1(x_i)+\frac{3}{4}\delta}^n$,
$w_n^{3,i}=1$ on $Q_i \setminus R_i$ for $n$ large, and
\begin{equation}
\label{mm3}
\limsup_n \sum_{Q_i \subseteq \Om} MM_n(w_n^{3,i})_{|\{v_n(t) > b_n^1\}} \le o(\sigma).
\end{equation}
Moreover there exists a cut-off function $\varphi_n^{3,i} \in H^1(Q_i)$ such that
$\varphi_n^{3,i}=0$ in a neighborhood of
$H_i^+ \setminus E_{a_2(x_i)-\delta}^n$ and of
$H_i^- \cap E_{a_1(x_i)+\delta}^n$, $\varphi_n^{3,i}=1$ on $Q_i \setminus R_i$ for
$n$ large,
$\supt{\nabla \varphi_n^{3,i}} \subseteq \{w_n^{3,i}=0\}$, and
\begin{equation}
\label{grad3}
\lim_n \eta_n \int_{Q_i \cap \{v_n(t)>b_n^1\}} |\nabla \varphi_n^{3,i}|^2\,dx=0.
\end{equation}
\end{lemma}

\begin{proof}
Let $\pi_i^\pm$ be the planes which contain $H_i^\pm$, and for $x \in \Om'$, let
$\pi_i^\pm x$ be its projection on $\pi_i^\pm$.
Let us now consider $(u_n(t))_{|H_i^+}$: we pose
$$
\psi^{i,+}_n(y):=\frac{4}{\delta}
\left( u_n(y)-a_2(x_i)+\frac{3}{4}\delta \right)^+ \wedge 1
$$
Note that $\psi^{i,+}_n$ is equal to zero on
$H_i^+ \setminus E^n_{a_2(x_i)-\frac{3}{4}\delta}$
and so on $\{x \in H_i^+\,:\, u_n(t)(x)= \gamma_n^i\}$ where $\gamma_n^i$ is defined
as in Lemma \ref{wni2}.
Moreover, $\psi_n^{i,+}=1$ on $H_i^+ \cap E_{a_2(x_i)-\frac{\delta}{2}}^n$.
If $d \in ]0,1]$
is such that $\eta_n+(b_n^1)^2 \ge d^2$, by (\ref{gradh+}) we have
\begin{equation}
\label{gradpsibound}
\int_{H_i^+ \cap \{v_n(t) > b_n^1\}} |\nabla \psi^{i,+}_n|^2 \,d\hs^{N-1} \le \frac{16}{\delta^2}
\int_{H_i^+ \cap \{v_n(t) > b_n^1\}} |\nabla u_n|^2 \,d\hs^{N-1} \le
\frac{16K_n}{\delta^2 d^2}.
\end{equation}
Let us define
$$
\tilde{\psi}^{i,+}_n(y):=\frac{4}{\delta}(u_n(y)-a_2(x_i)+\delta)^+ \wedge 1
$$
which is null on $H_i^+ \setminus E_{a_2(x_i) -\delta}^n$.
\par
In a similar way we construct
$\psi_n^{i,-}$ and $\tilde{\psi}^{i,-}_n$ on $H_i^-$ which are null on
$H_i^- \cap E^n_{a_1(x_i)+\frac{3}{4}\delta}$ and on
$H_i^- \cap E^n_{a_1(x_i)+\delta}$ respectively.
Let us pose
\begin{equation*}
w_n^{3,i,\pm}(x):=
\left[
\psi_n^{i, \pm}(\pi_i^\pm x) +\frac{1}{\varepsilon_n}(d_{H_i^\pm}(x)-l_n^i)^+
\right] \wedge 1
\end{equation*}
with $\frac{l_n^i}{\varepsilon_n} \to 0$ and $\frac{\eta_n}{(l_n^i)^2} \to 0$. This
is possible since $\eta_n << \varepsilon_n$.
Let
$A_n^i:= (H_i^+ \setminus E_{a_2(x_i)-\frac{\delta}{2}}^n)
\times ]-\varepsilon_n-l_n^i,\varepsilon_n+l_n^i[
\cap \{v_n(t) > b^1_n\}$.
Then we have by definition of $\psi_n^{i, \pm}$, by (\ref{h+est}),
(\ref{gradpsibound}) and the fact that $K_n \varepsilon_n$ is bounded in $n$
\begin{multline*}
\limsup_n MM_n(w_n^{3,i,+})_{|\{v_n(t) > b^1_n\}}= \\
=\limsup_n
\left\{
\frac{\varepsilon_n}{2}
\int_{Q_i \cap \{v_n(t) > b^1_n\}} |\nabla w_n^{3,i,+}|^2 \, dx
+\frac{1}{2\varepsilon_n}
\int_{Q_i \cap \{v_n(t) > b^1_n\}} (1-w_n^{3,i,+})^2 \,dx
\right\} \le \\
\le \limsup_n
\left\{
\frac{\varepsilon_n}{2}
\int_{A_n^i}
\left(|\nabla \psi_n(\pi_i^+ x)|^2 +\frac{1}{\varepsilon_n^2} \right) \,dx+ \right.
\\ \left.
+\frac{1}{2\varepsilon_n}
\left( 2\hs^{N-1}(H_i^+ \setminus E_{a_2(x_i)-\frac{\delta}{2}}^n)
(\varepsilon_n+l_n^i)
\right) \right\}
\end{multline*}
so that we get
\begin{multline*}
\limsup_n MM_n(w_n^{3,i,+})_{|\{v_n(t) > b^1_n\}}= \\
\le \limsup_n
\left\{
\frac{\varepsilon_n}{2} \frac{16K_n}{\delta^2d^2} 2(\varepsilon_n+l_n^i)
+\frac{\varepsilon_n}{2} \frac{2}{\varepsilon_n^2} (\varepsilon_n+l_n^i)
\hs^{N-1}(H_i^+ \setminus E_{a_2(x_i)-\frac{\delta}{2}}^n) + \right. \\
\left.
+\frac{\varepsilon_n+l_n^i}{\varepsilon_n}
\hs^{N-1}(H_i^+ \setminus E_{a_2(x_i)-\frac{\delta}{2}}^n)
\right\} \le
\\ \le 2\limsup_n
\hs^{N-1}(H_i^+ \setminus E_{a_2(x_i)-\frac{\delta}{2}}^n) \le 4\sigma r_i^{N-1}.
\end{multline*}
Similar calculations hold for $w_n^{3,i,-}$. Let us pose $w_n^{3,i}:=
w_n^{3,i,+} \wedge w_n^{3,i,-}$.
Then $0 \le w_n^{3,i} \le 1$, $w_n^{3,i}=0$ in a neighborhood of
$H_i^+ \setminus E_{a_2(x_i)-\frac{3}{4}\delta}^n$ and of
$H_i^- \cap E_{a_1(x_i)+\frac{3}{4}\delta}^n$,
$w_n^{3,i}=1$ on $Q_i \setminus R_i$ for $n$ large, and we have that
\begin{equation*}
\limsup_n \sum_{Q_i \subseteq \Om}
MM_n(w_n^{3,i})_{|\{v_n(t) > b_n^1\}} \le o(\sigma),
\end{equation*}
which prove the first part of the lemma.
\par
We define
\begin{equation*}
\varphi_n^{3,i}(x):=
\left[
\tilde{\psi}^{i,+}_n(\pi_i^+ x) +\frac{1}{l_n^i}
\left( d_{H_i^+}(x)-\frac{l_n^i}{2} \right)^+
\right]
\wedge
\left[
\tilde{\psi}^{i,-}_n(\pi_i^- x) +\frac{1}{l_n^i}
\left( d_{H_i^-}(x)-\frac{l_n^i}{2} \right)^+
\right]
\wedge 1.
\end{equation*}
The previous calculations prove that
\begin{equation*}
\lim_n \eta_n \int_{Q_i \cap \{v_n(t)>b^1_n\}} |\nabla \varphi_n^{3,i}|^2 \,dx = 0
\end{equation*}
since $\displaystyle \frac{\eta_n}{(l_n^i)^2} \to 0$. Moreover
$\varphi_n^{3,i}=1$ on $Q_i \setminus R_i$ for $n$ large.
\end{proof}

\begin{lemma}
\label{w4}
Suppose that $Q_i \subseteq \Om$; then there exists $w_n^{4,i} \in H^1(\Om')$ such
that $0 \le w_n^{4,i} \le 1$, $w_n^{4,i}=0$ in a neighborhood of $V_i$,
$w_n^{4,i}=1$ on $\Om_D$ for $n$ large
and
\begin{equation}
\label{mm4}
\limsup_n \sum_{Q_i \subseteq \Om}
MM_n(w_n^{4,i})\le o(\sigma).
\end{equation}
Moreover there exists a cut-off function $\varphi_n^{4,i}$ such that
$\varphi_n^{4,i}=0$ in a
neighborhood of $V_i$, $\varphi_n^{4,i}=1$ on $\Om_D$ for $n$ large,
$\supt{\nabla \varphi_n^{4,i}} \subseteq \{w_n^{4,i}=0\}$,
and
\begin{equation}
\label{grad4}
\lim_n \eta_n \int_{\Om'} |\nabla \varphi_n^{4,i}|^2 \,dx =0.
\end{equation}
\end{lemma}

\begin{proof}
Let us pose
$$
w_n^{4,i}(x):=
\frac{1}{\varepsilon_n} (d_{V_i}(x)-l^i_n)^+ \wedge 1,
$$
and
$$
\varphi_n^{4,i}(x):=
\frac{1}{l^i_n} \left( d_{V_i}(x)-\frac{l^i_n}{2} \right)^+ \wedge 1
$$
where $\frac{l_n^i}{\varepsilon_n} \to 0$ and $\frac{\eta_n}{(l_n^i)^2} \to 0$.
We have immediately (since $\sum_{Q_i \subseteq \Om} \hs^{N-1}(V_i) \le o(\sigma)$)
\begin{equation*}
\limsup_n \sum_{Q_i \subseteq \Om}
MM_n(w_n^{4,i})\le o(\sigma)
\end{equation*}
while
\begin{equation*}
\lim_n \eta_n \int_{\Om'} |\nabla \varphi_n^{4,i}|^2 \,dx = 0
\end{equation*}
since $\frac{\eta_n}{(l^i_n)^2} \to 0$. For $n$ large enough,
$w_n^{4,i}=1$, $\varphi_n^{4,i}=1$ on $\Om_D$ and the proof is complete.
\end{proof}

We recall that $z=g_h(t)$ in a neighborhood $\mathcal V$ of
$\partial \Om \setminus \cup Q_i$.

\begin{lemma}
\label{w4bdry+}
Let $Q_i \cap \partial_D \Om \not= \emptyset$ with $Q_i^+ \setminus R_i \subseteq \Om$.
Then $E_{a_1(x_i)+\frac{\delta}{2}}^n \cap Q_i \subseteq \Om$ for all $n$, and
there exists $w_n^{b,i,+} \in H^1(\Om')$ with  $0 \le w_n^{b,i,+} \le 1$,
$w_n^{b,i,+}=1$ on $\Om_D$, $w_n^{b,i,+}=0$ in a neighborhood of
$$
V_i^{n,+}:=\left[ V_i \cap E_{a_1(x_i)+ \delta}^n \right]
\cup \left[ (V_i \cap Q_i^+) \setminus {\mathcal V} \right],
$$
and such that
\begin{equation}
\label{mm4bdry+}
\limsup_n \sum_{Q_i \cap \partial_D \Om \not= \emptyset}
MM_n(w_n^{b,i,+})_{|\{v_n(t) > b_n^1\}} \le o(\sigma).
\end{equation}
Moreover there exists a cut-off function $\varphi_n^{b,i,+}$
such that $\varphi_n^{b,i,+}=1$ on $\Om_D$, $\varphi_n^{b,i,+}=0$
in a neighborhood of $V_i^{n,+}$,
$\supt{\nabla \varphi_n^{b,i,+}} \subseteq \{w_n^{b,i,+}=0\}$, and
\begin{equation}
\label{grad4bdry+}
\lim_n \eta_n \int_{\Om' \cap \{v_n(t)>b_n^1\}} |\nabla \varphi_n^{b,i,+}|^2 \,dx =0
\end{equation}
\end{lemma}

\begin{proof}
Note that by construction, $E_{a_1(x_i)+\frac{\delta}{2}}^n \cap Q_i \subseteq \Om$ since
$u_n(t)$ is continuous and $u_n(t)=g_h(t)$ on $\Om_D$.
It is now sufficient to operate as in Lemma \ref{wn3} and in Lemma \ref{w4}.
In fact, in view of \eqref{boundvu}, we may construct
$\tilde{w}_n^{b,i,+} \in H^1(\Om')$ such that $0 \le \tilde{w}_n^{b,i,+} \le 1$,
$\tilde{w}_n^{b,i,+}=0$ in a neighborhood of
$V_i \cap E_{a_1(x_i)+\delta}^n$, $\tilde{w}_n^{b,i,+}=1$ on $\Om_D$ and on
$V_i \setminus E_{a_1(x_i)+\frac{\delta}{2}}^n$, and such that
$\limsup_n MM_n(\tilde{w}_n^{b,i,+})_{|\{v_n(t) > b_n^1\}} \le o(\sigma)r_i^{N-1}$.
Referring to
$(V_i \cap Q_i^+) \setminus {\mathcal V}$, we can reason as in Lemma \ref{w4} getting
$\overline{w}_n^{b,i,+}$, such that $0 \le \overline{w}_n^{b,i,+} \le 1$,
$\overline{w}_n^{b,i,+}=0$ in a neighborhood of
$(V_i \cap Q_i^+) \setminus {\mathcal V}$, $\overline{w}_n^{b,i,+}=1$ on $\Om_D$,
and such that $\limsup_n MM_n(\overline{w}_n^{b,i,+})\le o(\sigma)r_i^{N-1}$.
\par
Posing $w_n^{b,i,+}:=\tilde{w}_n^{b,i,+} \wedge \overline{w}_n^{b,i,+}$, we get the
first part of the thesis. Similarly, we may construct $\varphi_n^{b,i,+}$ which
satisfies (\ref{grad4bdry+}).
\end{proof}

In a similar way we can prove the following lemma.

\begin{lemma}
\label{w4bdry-}
Let $Q_i \cap \partial_D \Om \not= \emptyset$ with
$Q_i^- \setminus R_i \subseteq \Om$.
Then $Q_i \setminus E_{a_2(x_i)-\frac{\delta}{2}}^n \subseteq \Om$ for all $n$, and
there exists $w_n^{b,i,-} \in H^1(\Om')$ with $0 \le w_n^{b,i,-} \le 1$,
$w_n^{b,i,-}=1$ on $\Om_D$, $w_n^{b,i,-}=0$ in a neighborhood of
$$
V_i^{n,-}:=\left[ V_i \setminus E_{a_2(x_i)- \delta}^n \right]
\cup \left[ (V_i \cap Q_i^-) \setminus {\mathcal V} \right],
$$
and such that
\begin{equation}
\label{mm4bdry-}
\limsup_n \sum_{Q_i \cap \partial_D \Om \not= \emptyset}
MM_n(w_n^{b,i,-})_{|\{v_n(t) > b_n^1\}} \le o(\sigma).
\end{equation}
Moreover there exists a cut-off function $\varphi_n^{b,i,-}$
such that $\varphi_n^{b,i,-}=1$ on $\Om_D$, $\varphi_n^{b,i,-}=0$
in a neighborhood of $V_i^{n,-}$,
$\supt{\nabla \varphi_n^{b,i,-}} \subseteq \{w_n^{b,i,-}=0\}$, and
\begin{equation}
\label{grad4bdry-}
\lim_n \eta_n \int_{\Om' \cap \{v_n(t)>b_n^1\}} |\nabla \varphi_n^{b,i,-}|^2 \,dx =0
\end{equation}
\end{lemma}

We can now prove Lemma \ref{minlem}.

\begin{proof}[Proof of Lemma \ref{minlem}]
We employ the notation of the preceding lemmas.
Following \cite[Theorem 2.1]{FL}, for each $i$ let us define $z_i^+$ on
$Q^+_i \cup R_i$ to be equal to $z$ on $Q_i^+ \setminus R_i$ and to the
symmetrization of $z$ with respect to $H_i(\sigma)$ on $R_i$.
Similarly we define $z_i^-$.
\par
For each $Q_i \subseteq \Om$, let us pose $z_n^i$ to be
equal to $z_i^+$ on $(Q_i^+ \setminus \tilde{R}_i) \cup
(E^n_{\gamma_n^i} \cap \tilde{R}_i)$, and to $z_i^-$ in the rest of $Q_i$.
\par
If $Q_i \cap \partial_D \Om \not= \emptyset$ with $Q_i^+ \setminus R_i \subseteq \Om$,
by Lemma \ref{wni2} and Lemma \ref{w4bdry+} we have
$E_{\gamma_n^i-\tau_n^i}^n \cap Q_i \subseteq Q_i^+$ for all $n$,
and its closure does not intersect $\partial \Om$.
We define $z_n^i$ to be equal to $z_i^+$ on
$(Q_i^+ \setminus \tilde{R}_i) \cup (E^n_{\gamma_n^i} \cap \tilde{R}_i)$, and
to $g_h(t)$ in the rest of $Q_i$. If $Q_i^- \setminus R_i \subseteq \Om$, by Lemma
\ref{wni2} and Lemma \ref{w4bdry-} we have
$Q_i \setminus E_{\gamma_n^i+\tau_n^i}^n  \subseteq \Om$, and
its closure does not intersect $\partial \Om$. We define $z_n^i$ to be equal to $z_i^-$
on $(Q_i^- \setminus \tilde{R}_i) \cup (\tilde{R}_i \setminus
E_{\gamma_n^i}^n)$, and to $g_h(t)$ in the rest of $Q_i$.
\par
Let us now define $\tilde{z}_n$ to be equal to $z$ outside
$\bigcup_{i=1}^k R_i$, and to $z_n^i$ inside each $R_i$.
We have $\tilde{z}_n=g_h(t)$ on $\Om_D$.
Note that if $Q_i \subseteq \Om$,
$H_i^+ \setminus E^n_{\gamma_n^i}$,
$H_i^- \cap E^n_{\gamma_n^i}$, $V_i^\pm$, and $\partial^*E^n_{\gamma_n^i} \cap Q_i$
could be contained in $S_{\tilde{z}_n}$.
Similarly, if $Q_i \cap \partial \Om \not= \emptyset$ and
$Q_i^+ \setminus R^i \subseteq \Om$ (the other case being similar), then
$H^+_i \setminus E_{\gamma_n^i}^n$, $V_i^{n,\pm}$ and
$\partial^*E^n_{\gamma_n^i} \cap Q_i$ could be contained in $S_{\tilde{z}_n}$.
\par
By assumption on $U$, we have that
\begin{equation}
\label{zapprox}
||\tilde{z}_n-z||_{L^2(\Om')}+||\nabla \tilde{z}_n- \nabla z||_{L^2(\Om';\R^N)} \le
o(\sigma);
\end{equation}
moreover, besides the possible jumps previously individuated,
$\tilde{z}_n$ has in $R_i$ polyhedral jumps which are a reflected version of the
polyhedral jumps of $z$ in $Q_i$. By assumption on $z$, we conclude that the union of
these polyhedral sets $P_i(S_z)$ has $\hs^{N-1}$ measure which is of the order of
$\sigma$ that is $\hs^{N-1}(P(S_z)) \le o(\sigma)$ where $P(S_z):= \bigcup_{i=1}^k P_i(S_z)$.
\par
Let $\tilde{w}_n$ be optimal for the Ambrosio-Tortorelli
approximation of $\left[ S_z \setminus (\bigcup Q_i) \right] \cup P(S_z)$ (as we
can find for example in \cite[Lemma 3.3]{Fo}),
that is $\tilde{w}_n$ is null in a neighborhood of
$\left[ S_z \setminus (\bigcup Q_i) \right] \cup P(S_z)$ and
\begin{align}
\label{mmreflected}
\limsup_n & MM_n(\tilde{w}_n) \le
\hs^{N-1}(S_z \setminus (\cup Q_i) \cup P(S_z)) \le \\
\nonumber
&\le \hs^{N-1}(S_z \setminus S_{u(t)}) +o(\sigma).
\end{align}
As in \cite{Fo}, let $\tilde{\varphi}_n$ be a cut-off function associated to
$\tilde{w}_n$, such that
\begin{equation}
\label{gradreflected}
\lim_n \eta_n \int_{\Om'} |\nabla \tilde{\varphi}_n|^2 \,dx =0.
\end{equation}
Let us pose for all $Q_i \subseteq \Om$
\begin{equation*}
w^i_n:=
\left\{
\begin{array}{l}
\min \{\tilde{w}_n, w_n^{2,i}, w_n^{3,i}, w_n^{4,i}\} \\ \\
\min \{\tilde{w}_n, w_n^{3,i}, w_n^{4,i}\} \\ \\
\min \{\tilde{w}_n, w_n^{4,i}\}
\end{array}
\begin{array}{l}
{\rm in}\;\tilde{R}_i \\ \\
{\rm in}\; R_i \setminus \tilde{R}_i \\ \\
{\rm outside}\; R_i,
\end{array}
\right.
\end{equation*}
and
\begin{equation*}
\varphi^i_n:=
\left\{
\begin{array}{l}
\min \{\tilde{\varphi}_n, \varphi_n^{2,i}, \varphi_n^{3,i}, \varphi_n^{4,i}\} \\ \\
\min \{\tilde{\varphi}_n, \varphi_n^{3,i}, \varphi_n^{4,i}\} \\ \\
\min \{\tilde{\varphi}_n, \varphi_n^{4,i}\}
\end{array}
\begin{array}{l}
{\rm in}\;\tilde{R}_i \\ \\
{\rm in}\; R_i \setminus \tilde{R}_i \\ \\
{\rm outside}\; R_i.
\end{array}
\right.
\end{equation*}
For all $Q_i$ such that $Q_i \cap \partial_D \Om \not=\emptyset$ with
$Q_i^+ \setminus R_i \subseteq \Om$, let us pose
\begin{equation*}
w^i_n:=
\left\{
\begin{array}{l}
\min \{\tilde{w}_n, w_n^{2,i}, w_n^{3,i,+}, w_n^{b,i,+}\} \\ \\
\min \{w_n^{2,i}, w_n^{b,i,+}\} \\ \\
\min \{\tilde{w}_n, w_n^{3,i,+}, w_n^{b,i,+}\} \\ \\
1 \\ \\
\min \{\tilde{w}_n, w_n^{b,i,+}\}
\end{array}
\begin{array}{l}
{\rm in}\;\tilde{R}_i \cap E_{\gamma_n^i}^n \\ \\
{\rm in}\;(\tilde{R}_i \setminus E_{\gamma_n^i}^n) \cup Q_i^- \\ \\
{\rm in}\;R_i \setminus (\tilde{R}_i \cup Q_i^-) \\ \\
{\rm in}\;\Om_D \\ \\
{\rm otherwise}
\end{array}
\right.
\end{equation*}
and
\begin{equation*}
\varphi^i_n:=
\left\{
\begin{array}{l}
\min \{\tilde{\varphi}_n, \varphi_n^{2,i}, \varphi_n^{3,i,+}, \varphi_n^{b,i,+}\} \\ \\
\min \{\varphi_n^{2,i}, \varphi_n^{b,i,+}\} \\ \\
\min \{\tilde{\varphi}_n, \varphi_n^{3,i,+}, \varphi_n^{b,i,+}\} \\ \\
1 \\ \\
\min \{\tilde{\varphi}_n, \varphi_n^{b,i,+}\}
\end{array}
\begin{array}{l}
{\rm in}\;\tilde{R}_i \cap E_{\gamma_n^i}^n \\ \\
{\rm in}\;(\tilde{R}_i \setminus E_{\gamma_n^i}^n) \cup Q_i^- \\ \\
{\rm in}\;R_i \setminus (\tilde{R}_i \cup Q_i^-) \\ \\
{\rm in}\;\Om_D \\ \\
{\rm otherwise}
\end{array}
\right.
\end{equation*}
Similarly we reason for the case $Q_i^- \setminus R_i \subseteq \Om$.
By construction, for all $i=1, \ldots,k$ we have that $w_n^i, \varphi_n^i \in H^1(\Om')$,
$0 \le w_n^i,\varphi_n^i \le 1$ and $w_n^i,\varphi_n^i=1$ on $\Om_D$ for $n$ large.
\par
Note that by Lemmas \ref{wni2}, \ref{wn3}, \ref{w4},
\ref{w4bdry+} and \ref{w4bdry-}, and by (\ref{mmreflected}) and (\ref{gradreflected}),
we have that
\begin{equation}
\label{mmnq}
\limsup_n \sum_{i=1}^k MM_n(w_n^i)_{|\{v_n(t) > b_n^1\}} \le
\hs^{N-1}(S_z \setminus S_{u(t)})+o(\sigma),
\end{equation}
and
\begin{equation}
\label{gradnq}
\lim_n \eta_n \sum_{i=1}^k
\int_{\Om' \cap \{v_n(t)>b_n^1\}} |\nabla \varphi_n^i(x)|^2 \,dx =0.
\end{equation}
We are now in a position to conclude the proof. We pose
$$
v_n:= \min \{w_n, w_n^i, i=1, \ldots, k\}, \quad
\varphi_n:= \min\{\varphi_n, \varphi_n^i, i=1, \ldots, k\}.
$$
Note that $\varphi_n=0$ in a neighborhood of $S_{\tilde{z}_n}$, and $\varphi_n=1$ on
$\Om_D$ for $n$ large. Moreover $0 \le v_n \le w_n \le v_n(t)$ in $\Om'$
and $v_n=1$ on $\Om_D$. Let $z_n:= \varphi_n \tilde{z}_n$; we have
$z_n \in H^1(\Om')$ with $z_n=g_h(t)$ on $\Om_D$. By (\ref{minimprop}), we have that
$$
\Atn{u_n(t)}{v_n(t)} \le \Atn{z_n}{v_n},
$$
and so
$$
\elln{u_n(t)}{v_n(t)} \le
\int_{\Om'} (\eta_n +v_n^2) |\nabla(\varphi_n \tilde{z}_n|^2 \,dx
+ MM_n(v_n)-MM_n(v_n(t)).
$$
We may write
\begin{align*}
\int_{\Om'} & (\eta_n +v_n(t)^2) |\nabla u_n(t)|^2 \,dx \le \\
&\le \int_{\Om'} (\eta_n+1) |\nabla \tilde{z}_n|^2 \,dx +
\int_{\Om'}
(\eta_n +v_n^2) (2\nabla \varphi_n \nabla \tilde{z}_n+
\tilde{z}_n |\nabla \varphi_n|^2) \,dx+ \\
& +MM_n(w_n)-MM_n(v_n(t)) +\sum_{i=1}^k MM_n(w_n^i)_{|\{v_n(t) > b_n^1\}}.
\end{align*}
Taking into account \eqref{zapprox}, \eqref{grad1}, \eqref{gradnq}, \eqref{mm1}, and \eqref{mmnq},
we have that passing to the limit
\begin{align*}
\int_{\Om'} |\nabla u|^2 \,dx \le
\int_{\Om'} &|\nabla z|^2 \,dx +\hs^{N-1}(S_z \setminus S_{u(t)}) +\\
&+\frac{2Ck}{(k-1)^2}+\frac{C}{(k-1)(1-b)^2}+\frac{Cb}{(1-b)^2}
+o(\sigma),
\end{align*}
so that, letting $\sigma \to 0$ and then $b \to 0$, $k \to \infty$ (which is permitted
choosing appropriately $j_2$ and $j_3$), we obtain the thesis.
\end{proof}

\vskip10pt
We can now pass to the proof of Theorem \ref{minthm}.
Given $0=t_1\le t_2 \le \ldots \le t_k=t$, it is sufficient to prove that
\begin{equation}
\label{mineqk}
\int_{\Om'} |\nabla u(t)|^2 \,dx \le
\int_{\Om'} |\nabla z|^2 \,dx +
\hs^{N-1} \left( S_z \setminus \left( \bigcup_{i=1}^k S_{u(t_i)} \right) \right).
\end{equation}
Passing to the sup on $t_1, \ldots, t_k$, we deduce in fact the thesis.
We obtain (\ref{mineqk}) using the same arguments of Lemma \ref{minlem}; defining
$$
J_j:=
\left\{ x \in \bigcup_{m=1, \ldots, k}
\left( \bigcup_{a_1,a_2 \in A^k}
\left[ \partial^* E_{a_1}^k \cap \partial^* E_{a_2}^k \right]
\right) \,:
\min_{l=1, \ldots, k} [u_l(x)] >\frac{1}{j} \right\},
$$
where $E_a^k$ and $A_k$ are defined as the corresponding sets for $u(t)$,
following \cite{FL}, we cover $J_j$ in such a way that for all $x_i \in J_j$ there
exists $l$ with $x_i \in S_{u(t_l)}$ and
$$
\hs^{N-1} \left( \left[ \bigcup_{r=1, \ldots, k} S_{u(t_r)} \setminus S_{u(t_l)} \right]
\cap Q_i \right) \le \sigma r_i^{N-1}.
$$
So in each $Q_i$ there exists $u(t_l)$ such that $\bigcup_{r=1}^k S_{u(t_r)} \cap Q_i$ is
essentially (with respect to the measure $\hs^{N-1}$) $S_{u(t_l)} \cap Q_i$.
Recalling that $v_n(t) \le v_n(t_l)$ for all $l=1, \ldots,k$, we have
$$
\ellnp{u_n(t_l)}{v_n(t)} \le
\ellnp{u_n(t_l)}{v_n(t_l)} \le C,
$$
and so it is readily seen that the arguments of Lemma \ref{minlem} can be adapted to prove
(\ref{mineqk}).

\section{A final remark}

The previous results can be extended to recover the case of non isotropic surface
energies, i.e., energies of the form
\begin{equation}
\label{nonisenergy}
\int_\Om |\nabla u|^2 \,dx+ \int_{\Gamma} \varphi(\nu_x) \,d\hs^{N-1}(x)
\end{equation}
where $\nu_x$ is the normal to $\Gamma$ at $x$, and $\varphi$ is a norm on $\R^N$.
In fact all the previous arguments are based on Theorem \ref{at} concerning the
elliptic approximation and on Theorem \ref{piecedensity} about the density of
piecewise smooth functions with respect to the elastic energy.
An elliptic approximation of Ambrosio-Tortorelli type of (\ref{nonisenergy})
has been proved in \cite{Fo}, while a density result of piecewise smooth functions
with respect to non-isotropic surface energies has been proved in \cite{CT}.
We conclude that all the previous theorems can be modified in order to treat the
more general energy (\ref{nonisenergy}).

\bigskip
\bigskip
\centerline{ACKNOWLEDGMENTS}
\bigskip\noindent
The author wishes to thank Gianni Dal Maso for having proposed him the problem,
and for many helpful and interesting discussions.

\end{document}